\documentclass[10pt,twoside]{article}

\usepackage[centertags]{amsmath}
\usepackage{amsfonts}
\usepackage{amssymb}
\usepackage{amsthm}
\usepackage{epsfig}
\usepackage{color}
\allowdisplaybreaks[4]

\newcommand{\COM}[1]{}
\definecolor{c20}{rgb}{0.,0.7,0.}
\definecolor{c30}{rgb}{0.,0.,1.}
\definecolor{c40}{rgb}{1,0.1,0.7}
\definecolor{c50}{rgb}{1,0,0}
\definecolor{c60}{rgb}{1,0.9,0.1}

\def\Tan#1{\textcolor{c20}{#1}}
\def\Tan#1{#1}

\def\cL#1{\textcolor{c50}{#1}}
 \def\cL#1{#1}

\def\ccL#1{\textcolor{c50}{#1}}
 \def\ccL#1{#1}

\newcommand{\kb}[1]{\boldsymbol{#1}}
\newcommand{\vk}[1]{\kb{#1}}

\newcommand{\ve}{\varepsilon}

\newcommand{\abs}[1]{\left\lvert #1 \right\rvert}
\newcommand{\Abs}[1]{ \biggl \lvert #1 \biggr \rvert}

\newcommand{\E}[1]{\mathbb{E}\left\{#1\right\}}

\newcommand{\pk}[1]{\mathbb{P} \left \{#1 \right \} }

\newcommand{\R}{\mathbb{R}}

\newcommand{\N}{\mathbb{N}}

\newcommand{\inn}{\in \N}

\newcommand{\limit}[1]{\lim_{#1 \to   \infty}}

\newcommand{\BQN}{\begin{eqnarray}}
\newcommand{\EQN}{\end{eqnarray}}
\newcommand{\BQNY}{\begin{eqnarray*}}
\newcommand{\EQNY}{\end{eqnarray*}}

\newcommand{\BS}{\begin{sat}}
\newcommand{\ES}{\end{sat}}
\newcommand{\BT}{\begin{theo}}
\newcommand{\ET}{\end{theo}}
\newcommand{\BL}{\begin{lem}}
\newcommand{\EL}{\end{lem}}
\newcommand{\BK}{\begin{korr}}
\newcommand{\EK}{\end{korr}}

\newcommand{\BD}{\begin{de}}
\newcommand{\ED}{\end{de}}

\newcommand{\BIT}{\begin{itemize}}
\newcommand{\EIT}{\end{itemize}}
\newcommand{\BDI}{\begin{description}}
\newcommand{\EDI}{\end{description}}

\newcommand{\BRM}{\begin{remarks}}
\newcommand{\ERM}{\end{remarks}}

\newcommand{\BEL}{\begin{lem}}
\newcommand{\EEL}{\end{lem}}

\newtheorem{theo}{Theorem}[section]
\newtheorem{sat}[theo]{Proposition}
\newtheorem{de}[theo]{Definition}
\newtheorem{lem}[theo]{Lemma}

\newtheorem{korr}[theo]{Corollary}
\newtheorem{remark}[theo]{Remark}
\newtheorem{remarks}[theo]{Remarks}

\newcommand{\nelem}[1]{{Lemma \ref{#1}}}

\newcommand{\netheo}[1]{{Theorem \ref{#1}}}

\newcommand{\prooflem}[1]{\textsc{\bf Proof of Lemma} \ref{#1}:}

\newcommand{\QED}{\hfill $\Box$}

\def\IF{\infty}

\newcommand{\expon}[1]{\exp\left(#1\right)}
\newcommand{\vnorm}[1]{\left\|#1\right\|}

\date{}

\def\H{\mathcal{H}}

\def\H{\mathcal{H}}

\def\X{\vk X}

\begin{document}
\baselineskip 15pt \setcounter{page}{1}
\title{\bf \Large  \ccL{On maxima of chi-processes over threshold dependent grids}
\thanks{Research supported by the Swiss National Science Foundation grant 200021-140633/1, RARE
-318984 (an FP7 Marie Curie IRSES Fellowship), Zhongquan  Tan also acknowledges National Science Foundation of China
(No. 11326175) and Natural Science Foundation of Zhejiang Province of China (No. LQ14A010012).} }
\author{{\small Chengxiu Ling$^a$, \quad Zhongquan  Tan$^{a,b}$ }\footnote{ E-mail address:  tzq728@163.com }\\
{\small\it $^a$Department of Actuarial Science, Faculty of Business and Economics (HEC Lausanne), University of Lausanne,}\\
{\small\it UNIL-Dorigny, 1015 Lausanne, Switzerland}\\
{\small\it  $^b$College of Mathematics, Physics and Information Engineering, Jiaxing University, Jiaxing 314001, PR China;}\\
}
 \maketitle
 \baselineskip 15pt

\begin{quote}
{\bf Abstract:}\ \ \COM{In this paper the asymptotic relation between the maximum of a
continuous dependent chi-process
and the maximum of this process sampled at discrete time points is
studied. It is shown that, for the weakly dependent case, these two maxima are
asymptotically independent, dependent and coincide when the grid of the discrete time
points is a sparse grid, Pickands grid and dense grid, respectively, while for the strongly dependent case,  these two maxima are
asymptotically totally dependent if the grid of the discrete time
points is sufficiently dense, and asymptotically dependent if the the grid points are sparse or Pickands grids.}
In this paper, with motivation from \ccL{\cite{Pit2004}} 
and the considerable interest in stationary chi-processes, we derive asymptotic joint distributions of maxima of stationary strongly dependent chi-processes on a continuous time and an uniform grid on the real axis. Our findings extend those for Gaussian cases  and give three  involved dependence structures via the strongly dependence condition and the sparse, Pickands and dense grids.



\end{quote}

\section{Introduction}\label{sec1}
Consider a stationary chi-process $\{\chi_m(t), t\ge0\}$ with $m, m\inn$ degrees of freedom as follows
\BQNY
\chi_m (t)=\left(X_{1}^{2}(t)+\cdots+X _{m} ^{2}(t)\right)^{1/2}= \|\X(t)\|,\quad t\ge 0,
\EQNY
where $\X(t)=(X_{1}(t),\ldots,X _{m} (t))$ is a vector Gaussian process which components are independent copies of a standard (zero-mean and unit-variance) stationary Gaussian process  $\{X(t), t\ge 0\}$ with almost surely (a.s.) continuous sample paths and correlation function $r(t)=\E{X(0) X(t)}, t\ge0$. \\
In this paper, we are concerned with the dependence of extremes of the continuous time and discrete time of chi-processes.
Specifically, assuming that the process $ \ccL{\{}\chi_m(t), t\in[0,T]\}$ is observed at time $t\in \mathfrak{R}(\delta)=\{k\delta, k\in\N\}$ with frequency $\delta=\delta_T>0$, of interest is the asymptotic joint distributions of
$( M_{m}(T),  M_m(\delta, T) )$ as $T\to\IF$ (after normalization) with
\BQN
\label{def:chi}
M_{m}(T):= \sup_{t\in[0,T]} \chi_m (t), \quad M_m(\delta, T):=\sup_{t\in\mathfrak{R}(\delta)\cap [0,T]} \chi_m (t).
\EQN
  The impetus for this investigation comes from numerical simulations of high extremes of continuous time random processes, see e.g.,
 \cite{Husler2004,Pit2004,TanH2014} for Gaussian processes, \cite{HuslerP2004} for the storage process with fractional
Brownian motion, \cite{HashorvaT2014Different, tan2012piterbarg, TanHP2013} for  stationary vector Gaussian processes and standardized stationary Gaussian processes, and \cite{TurkmanA} for stationary processes. It is shown in the aforementioned contributions that the dependence between continuous time extremes and discrete time  extremes is determined strongly by the sampling frequency $\delta$ and the normalization constants, see also for related discussions \cite{BrigoM2000Finance,LedfordT2003diagnostics,RobinsonT2000,Scotto2003subsampleTS, TurkmanA} in the financial and time series literature. Another motivation is that since  the chi-processes appear naturally as limiting processes which have attracted considerable interest from both theoretical and practical fields, see e.g., \cite{HashorvaJ2014trend, Liu2014extremes, Hashorva2014inflatedchi, HashorvaJi2013locally} for deeply theoretical discussions involved in the continuous time extremes of various $\chi$-processes, and  \cite{Albin2003,Aue2009, JaruskovaP02011, Jaruvskova2015detecting} for statistics test applications concerning the maxima  over the chosen time points set and the continuous time intervals. Therefore, of crucial importance is to understand the underlying asymptotic behavior of the extremes for different grids.\\
The principle challenge for $\chi$-processes
increases significantly due to no counterpart of \emph{Berman\rq{}s Normal Comparison
Lemma} for chi-distributions. However, with the technical methodology from \cite{HashorvaJ2014trend,Lindgren1980,Lindgren1984,Lindgren1989,
Pit96,TanH2013randominterval, PiterbargS2004Chi, MittalY1975},
and assuming certain locally and long range dependence on the common Gaussian process $ X(\cdot)$, namely  (see for its extensional utilizations \cite{MittalY1975, leadbetter1983extremes, tanH2012})
\BQN\label{corrr}
 r(t)= 1 -  \abs{t}^{\alpha}+ o(\abs{t}^{\alpha}), \quad  t\to 0\ \ {\rm for\ some }\ \alpha \in (0,2]
\EQN
and
\begin{eqnarray}
\label{cond.Global}
\limit T r(T)\ln T= r\in [0,\infty),
\end{eqnarray}
we establish our findings in \netheo{T1} extending those for weakly dependent stationary Gaussian processes in \cite{Pit2004}, corresponding to $m=1$ and $r=0$ in \eqref{cond.Global} in our setting, represent asymptotically completely dependence, max-stable dependence and conditional independence according to the three different types of grids in the terminology of \cite{Pit2004}, namely the dense grid $\mathfrak{R}(\delta)$ with $\delta(T)=o((2\ln T)^{-1/\alpha}), T\to\IF$, the Pickands grid $\mathfrak{R}(\delta)$ with $\delta(T)=D(2\ln T)^{-1/\alpha}$ for some $D\in(0,\IF)$,  and the sparse grid $\mathfrak{R}(\delta)$ with $\lim_{T\to\IF} \delta(T) (2\ln T)^{1/\alpha}=\IF$ and $\delta(T) \le \delta_0$ for some $\delta_0>0$.\\
We note in passing that our methodology is different from that in \cite{TurkmanA} which is strongly based on the Albin\rq{}s methodology wherein the verification of technical Albin\rq{}s conditions requires in general a lot of efforts. Moreover, our theoretical results, which do not seem possible to be guessed, are of interest for simulation studies, and give to some extent certain recommendations how tight a simulation grid should be when high extremes are important in simulations of the chi-processes under consideration. \\

\COM{then for any fixed $T>0$ 
\begin{eqnarray}
\label{def:chi}
\pk{\sup_{t\in[0,T]}\chi_m(t) >u}=T\frac{2^{1-m/2}\H_{\alpha}}{\Gamma(m/2)}u^{2/\alpha+m-2}\expon{-\frac{u^2}2}(1+o(1)), \quad  u\rightarrow\infty,
\end{eqnarray}
where $\Gamma(\cdot)$ is the Euler Gamma function and $\H_\alpha $ is the well-known Pickands constant given by  
$$\mathcal{H}_{\alpha}=\lim_{\lambda\rightarrow\infty} \frac{ \mathbb{E}\exp\left(\max_{t\in[0,\lambda]} B^*_{\alpha/2}(t)
\right)}{\lambda} \in(0,\IF),$$
see for  extensive discussions \cite{Husler1999extremes,Harper3, Shao1996bounds,
Dieker2014asymptotic}. Here $B_{H}$ is a fractional Brownian motion (fBm) with Hurst index $H\in(0,1]$, i.e., a Gaussian process with stationary increments such that
$\E{B_{H}^{2}(t)}=|t|^{2H}, \E{B_H(t)}= -\abs{t}^{2H}$. The asymptotic properties of $M_{m}(T)$  have been extensively studied in the literature; see \cite{Aronowich1985behaviour,  ArendarczykD2012GaussianRI, Lindgren1989,Piterbarg1994nonstationary, Sharpe1978chiUpcrossing, tanH2012, TanW2014ChiRI,DebickiHJ2014RI} for
various results.

}

\COM{\cite{Pit2004}  first studied the asymptotic relation between
$M_{T}$ and the maximum of the discrete version $M_{T}^{\delta}=
\max\{X(k\delta), 0 \leq k\delta \leq T \}$ for some
$\delta=\delta(T)>0, k\in \mathbb{N}$, where $\mathbb{N}$ denotes the set of all natural numbers. Following \cite{Pit2004},
we consider uniform grids
$\mathfrak{R}=\mathfrak{R}(\delta)=\{k\delta: k\in \mathbb{N}\}$,
$\delta>0$. A grid is called sparse if $\delta$ is such that
$$\delta(2\ln T)^{1/\alpha}\rightarrow D$$
with $D=\infty$. If $D\in(0, \infty)$, the grid is a Pickands grid, and if $D = 0$, the grid is dense.

For the stationary Gaussian processes, \cite{Pit2004}   first showed that the maximum
$M_{T}^{\delta}$ of discrete time points and the maximum $M_{T}$ of
the continuous time points can be asymptotically independent,
dependent or totally dependent if the grid is a sparse, a Pickands
or a dense grid, respectively. This type of results are called
Piterbarg's max-discretisation theorems in the literature, see e.g.,  Tan and Hashorva (2014).

Piterbarg's max-discretisation theorems have been
extended to more general Gaussian cases, see \cite{Husler2004} for the locally stationary
Gaussian case, \cite{HuslerP2004} for the storage process with fractional
Brownian motion, Tan and Hashorva (2014) for multivariate weakly and strongly dependent Gaussian case and Hashorva and Tan (2015) for different grids
Although the Piterbarg's max-discretisation theorems for Gaussian processes
have been  studied extensively under different conditions in the past, it is far from complete.
However, all of the above mentioned papers dealt with the Gaussian case.
Extending the above results to no-Gaussian case is also interesting, since most of reality can not
be modeled by Gaussian model.  Turkman (2012) considered this problem by adopting the model from Albin (1990), which
dealt with the extremes properties of no-Gaussian processes under some mixing condition like the well-known $D(u)$ condition
from \cite{leadbetter1983extremes}.
In this paper, we are interested in the similar problems for chi-processes. However Albin's model is very
general and it can include the weakly dependent chi-processes, so by applying the results of Turkman (2012), we can obtain the
Piterbarg's max-discretisation theorems for weakly dependent chi-processes. It seem that the strongly dependent case can not been derived from
the results of Turkman (2012) directly. The goal of this paper is to study the Piterbarg's max-discretisation theorems for strongly dependent chi-processes.}

The rest of the paper is organized as follows\ccL{.} Our main results are presented in the next section. All the proofs are relegated to Section \ref{sec3}  which is followed by an Appendix including some technical auxiliary results. 

\section{Main results}\label{sec2}This section is devoted to the asymptotic properties of $(M_m(T), M_m(\delta, T))$ given in \eqref{def:chi} for the three different types of grids $\delta=\delta(T)$ in the  terminology of \cite{Pit2004}. Before giving our main resul\ccL{t} 
 (see \netheo{T1} below), we shall first recall some asymptotic results of the considered chi-processes and introduce some notation concerning the Pickands type constants.  \\
As we know  from \cite{Piterbarg1994nonstationary} or Corollary 7.3 in \cite{Pit96} that, if  the correlation function $r(t)$ satisfies \eqref{corrr}  and in addition $r(t)<1$ for all $t\neq0$, then for any fixed $T>0$ 
\begin{eqnarray}
\label{def:chi0}
\pk{M_m(T) >u}=T\frac{2^{1-m/2}\H_{\alpha}}{\Gamma(m/2)}u^{2/\alpha+m-2}\expon{-\frac{u^2}2}(1+o(1)), \quad  u\rightarrow\infty,
\end{eqnarray}
where $\Gamma(\cdot)$ is the Euler Gamma function and $\H_\alpha\in(0,\IF) $ denotes  the Pickands constant,  
see \cite{Husler1999extremes,Harper3,Pit96, leadbetter1983extremes, Shao1996bounds,
Dieker2014asymptotic} for details and various discussions.
The asymptotic properties of $M_{m}(T)$  have been extensively studied in the literature; see \cite{ArendarczykD2012GaussianRI,Aronowich1985behaviour, DebickiHJ2014RI, Lindgren1989,Piterbarg1994nonstationary, Sharpe1978chiUpcrossing, tanH2012, TanW2014ChiRI} for
various results. Moreover, if additionally condition \eqref{cond.Global} holds for some $r\in[0,\IF)$, then the mixed Gumbel limit theorem holds as follows (see e.g., Theorem 3.1 in  \cite{TanH2013randominterval})
 \begin{eqnarray}
\label{eq TB}
\pk{a_{T}(M_m(T)-b_{T})\leq x}\to \E{\exp\left(-e^{-x-r+\sqrt{2r}\chi_{m}}\right)},\quad T\to\IF,
\end{eqnarray}
with  $\chi_m$  positive such that $\chi_m^2$ a chi-square random variable with $m$ degrees of freedom, and $a_{T}, b_{T}$  given by
\BQN\label{def:normalization.ab}
a_T=\sqrt{2 \ln T}, \quad b_T=a_T+ \frac{\ln \Bigl( 2^{1- m/2}(\Gamma(m/2))^{-1}\H_{\alpha}a_T^{2/\alpha+m-2}\Bigr)}{a_T}.
\EQN
\COM{Define a stationary  $\chi$-processes with $m\geq 2$ degrees of freedom as
\begin{eqnarray}
\label{def:chi}
\chi_m (t)=\vnorm{\X(t)}=\left(X_{1}^{2}(t)+\cdots+X _{m} ^{2}(t)\right)^{1/2},\ \ t\ge 0,
\end{eqnarray}
where $\X(\mathbf{t})=(X_{1}(t),\ldots,X _{m} (t))$ is a Gaussian vector process, the components of which are
independent copies of a standard stationary Gaussian process $\{X(t), t\ge 0\}$ with  correlation function $r(t)$.
For any grid $\mathfrak{R}(\delta)$, define
$$M_{m}(T):= \sup_{t\in[0,T]} \chi_m (t) \ \ \mbox{and}\ \  M_m(\delta, T):=\sup_{t\in\mathfrak{R}(\delta)\cap [0,T]} \chi_m (t).$$
By using different methods, the tail probability of $M_{m}(T)$
has been studied by many authors, see Sharpe (1978), Lindgren
(1980a,1980b, 1989), Aronowich and Adler (1985, 1986), Albin (1990),
Piterbarg (1994) and \cite{Pit96} for various results.
If $r(t)$ satisfies the condition (\ref{corrr}) by Corollary 7.3 \cite{Pit96}
for any fixed $T>0$ we have}
\COM{\begin{eqnarray}
\label{def:chi}
P\left(M_{m}(T) >u\right)=T\frac{2^{1-m/2}\H_{\alpha}}{\Gamma(m/2)}u^{2/\alpha+m-2}\exp(-u^2/2)(1+o(1)), \quad  u\rightarrow\infty,
\end{eqnarray}
where $\Gamma(\cdot)$ is the Euler gamma function. However, the
case when $T$ increases not too fast still holds, see  Theorem 7.2 of \cite{Pit96} for details.
The limit distribution of $M_{m}(T)$
has also been studied by many authors, see e.g.,, Konstant et al. (2004), \cite{PiterbargS2004Chi, StamatovicS2010,TanW2014ChiRI} for various results. A general limit distribution  of $M_{k}(T)$ is as follows.
Suppose that $r(t)$ satisfies the conditions (\ref{corrr}) and (\ref{cond.Global}),
by Tan and Hashorva (2013a), we have
\begin{eqnarray}
\label{eq TB}
P\{a_{T}(M_{m}(T)-b_{T})\leq x\}\to \mathbb{E}\exp\left(-e^{-x-r+\sqrt{2r}\chi_m }\right)
\end{eqnarray}
as $T\rightarrow\infty$, where Here $\chi_m$ is positive such that $\chi_m^2$ a chi-square random variable with $m$ degrees of freedom, and $a_{T}, b_{T}$ are given by
\BQN\label{def:normalization.ab}
a_T=\sqrt{2 \ln T}, \quad b_T=a_T+ a_T^{-1}\ln \Bigl( 2^{1- m/2}(\Gamma(m/2))^{-1}\H_{\alpha}a_T^{2/\alpha+m-2}\Bigr).\EQN
}
Next, we shall state our main result which \ccL{is} a type of Piterbarg's max-discretisation theorems  for chi-processes in terms of \cite{tan2012piterbarg}. To this end,  \ccL{} two Picaknds type constants (see \eqref{def: H.D} below) are needed. \\
Let $B^*_{\alpha/2}(t):=\sqrt 2 B_{\alpha/2}(t) -t^\alpha$ with $B_H(\cdot)$ a fractional Brownian motion (fBm) with Hurst index $H\in(0,1]$, and thus define a $m$-parameter fBm $B^*_{\boldsymbol{\alpha}/2}(\mathbf{t})=\sum_{i=1}^{m}B^*_{\alpha_{i}/2}(t_{i}), \boldsymbol{\alpha}=(\alpha_1, \ldots, \alpha_m)\in(0,2]^m, \mathbf t=(t_1, \ldots, t_m)\in[0,\IF)^m$ with mutually independent fBms $B^*_{\alpha_{i}/2}(\cdot), i\le m$, see e.g., \cite{Lin2001multifBm, Pit96} for related discussions on the $m$-parameter fBm $B^*_{\boldsymbol{\alpha}/2}(\cdot)$. We define thus, for any $D>0$ and $\boldsymbol {\alpha}_0=(\alpha, 2,\ldots,2)\ccL{\in(0,2]^m}$
\BQN\label{def: H.D}
\mathcal{H}_{D,\alpha}=\lim_{\lambda\rightarrow\infty}\frac{\mathbb{E}\exp\left(\max_{kD\in[0,\lambda], k\inn} B^*_{\alpha/2}(kD)\right)}{\lambda}, \quad \H_{D,\boldsymbol \alpha_0}^{x,y}:=\lim_{\lambda\rightarrow\infty}
\frac{\H_{D,\boldsymbol \alpha_0}^{x,y}(\lambda)}{\lambda^{m}}, \ x,y\in\R,
\EQN
which are finite and positive by Theorem 2 in \cite{Pit2004} and \nelem{L3.4}, respectively,  
here
\BQNY\H_{D,\boldsymbol{\alpha}}^{x,y}(\lambda)=\int_{-\infty}^{+\infty}e^{s}
\pk{\max_{\mathbf t\in [0,\lambda]^{m}} B^*_{\boldsymbol{\alpha}/2}(\mathbf t)>s+x,
\max_{\mathbf t\in [0,\lambda]^{m}\cap(\{kD, k\inn\}\times \R^{m-1})} B^*_{\boldsymbol{\alpha}/2}(\mathbf t)>s+y}ds,
\EQNY
and further 
 \BQN\label{def:normalization}
b_{\delta, T}=\left\{\begin{array}{ll}
a_T+ \frac{\ln \Bigl( 2^{1- m/2}(\Gamma(m/2))^{-1}\mathcal{H}_{D,\alpha}a_T^{2/\alpha+m-2}\Bigr)}{a_T}, & \mathfrak{R}(\delta){ \rm \ a \ Pickands\  grid };\\ 
a_T+ \frac{\ln \Bigl( 2^{1- m/2}(\Gamma(m/2))^{-1} \delta^{-1}a_T^{m-2}\Bigr)}{a_T}, & \mathfrak{R}(\delta){ \rm \ a \ sparse\  grid}.
\end{array}\right.
 \EQN

\BT\label{T1} Let $( M_{m}(T),  M_m(\delta, T) )$ be given as  in \eqref{def:chi}. Suppose that the correlation $r(\cdot)$  satisfies condition \eqref{corrr} and \eqref{cond.Global}, we have, with involved quantities given by \eqref{eq TB}--\eqref{def:normalization}, as $T\to\IF$ and $x, y\in\R$ \\
(a) For the sparse grid $\mathfrak R(\delta)$
\begin{eqnarray}
\label{eq2.1}
\pk{a_{T}\big(M_{m}(T)-b_{T}\big)\leq x,
           a_{T}\big( M_m(\delta, T)- b_{\delta, T}\big)\leq y}\to \mathbb{E}\exp\left(-\big(e^{-x}+e^{-y}\big)e^{-r+\sqrt{2r}\chi_{m}}\right).
\end{eqnarray}
(b) For the Pickands grid $\mathfrak{R}(\delta)=\mathfrak{R}(D(2\ln
T)^{-1/\alpha})$ with $D>0$
\begin{eqnarray}
\label{eq2.2}
\lefteqn{\pk{a_{T}\big(M_{m}(T)-b_{T}\big)\leq x,
           a_{T}\big( M_m(\delta, T)- b_{\delta, T}\big)\leq y}}\nonumber\\
&&\to \mathbb{E}\exp\left(-\big(e^{-x}
+e^{-y}-\pi^{(m-1)/2}\mathcal{H}_{D,\boldsymbol {\alpha}_0}^{\ln\mathcal{H}_{\alpha}+x, \ln
\mathcal{H}_{D,\alpha}+y}\big)e^{-r+\sqrt{2r}\chi_m }\right).
\end{eqnarray}
(c) For any dense grid $\mathfrak{R}(\delta)$
\begin{eqnarray}
\label{eq2.3}
&&\pk{a_{T}\big(M_{m}(T)-b_{T}\big)\leq x, a_{T}\big( M_m(\delta, T)- b_{T}\big)\leq y}
\to \mathbb{E}\exp\left(-e^{-\min(x,y)-r+\sqrt{2r}\chi_m }\right).
\end{eqnarray}
\ET


{\remark
(a) A straightforward application of \netheo{T1} (a) with $\delta(T)\equiv 1$ yields that
\begin{eqnarray*}
\label{eq2.4}
\pk{a_{T}\big(M_m(1, T)- b_{1,T}\big)\leq x}
\to \mathbb{E}\exp\left(-e^{-x-r+\sqrt{2r}\chi_m }\right), \quad x\in\R,
\end{eqnarray*}
which may have independent interest  in viewpoint of statistics applications,
see e.g., \cite{chareka2006test} for utilizations of the above limit  with $m=1$ concerning  test for additive outliers. \\
(b) From our results we  see that the joint convergence is determined by the choice of the grids and the normalization constants $a_T, b_T$ and $b_{\delta, T}$, which \ccL{is helpful} 
in  simulation studies and statistical applications, see related discussions for vector Gaussian processes in \cite{HashorvaT2014Different}. \\
(c) Clearly, the marginal distributions are the same, i.e., the mixed Gumbel distributions, and our results extend those for the Gaussian processes, see \cite{HashorvaT2014Different,Pit2004}. Moreover, the joint limit distribution for the Pickands grid is more involved due to the complication of the Pickands type constant $\H_{D, \boldsymbol \alpha_0}^{\ln \H_\alpha+x, \ln \H_{D,\alpha}+y}$, which calculation and simulation are open problems.\\ %
(d) It might be possible to allow $X_i\rq{}s$ to be dependent with condition \eqref{corrr} stated in a slightly general form such as $r_i(t)=1-C_i\abs{t}^{\alpha_i}(1+o(1)), \ t\to0$ as well. Results for extremes of chi-type processes for such generalizations can be found in \cite{Albin2003, Aue2009, Lindgren1989, Liu2014extremes}.\\
(e) It might be interesting to investigate the limit theorems for different grids as in \cite{HashorvaT2014Different}. Another possibility is to relax $r\in[0,\IF] $ in \eqref{cond.Global}; see e.g., \cite{MittalY1975, HashorvaT2014Different} for similar discussions. \\
(f) Following our arguments, it might be possible to consider the same problem for locally stationary chi-processes and cyclo-stationary chi-processes which are considered in  \cite{Husler2004, HashorvaJi2013locally} and \cite{Konstantinides2004gnedenko, tanH2012}, respectively.
}

\COM {\textbf{Theorem 2.2}. {\sl Let $\{X(t), t\geq0\}$ be a standard (zero-mean, unit-variance)
stationary Gaussian process with correlation functions $r(\cdot)$ satisfying (\ref{corrr}) and (\ref{cond.Global}). Then for any Pickands grid $\mathfrak{R}(\delta)=\mathfrak{R}(D(2\ln
T)^{-1/\alpha})$ with $D>0$,
\begin{eqnarray}
\label{eq2.2}
&&\pk{a_{T}\big(M_{m}(T)-b_{T}\big)\leq x,
           a_{T}\big( M_m(\delta, T)- b_{\delta, T}\big)\leq y}\nonumber\\
&&\ \ \ \ \ \to \mathbb{E}\exp\left(-\big(e^{-x-r+\sqrt{2r}\chi_m }
+e^{-y-r+\sqrt{2r}\chi_m }-\pi^{(m-1)/2}\mathcal{H}_{D,\boldsymbol {\alpha}_0}^{\ln\mathcal{H}_{\alpha}+x, \ln
\mathcal{H}_{D,\alpha}+y}e^{-r+\sqrt{2r}\chi_m }\big)\right)
\end{eqnarray}
as $T\rightarrow\infty,$ where $\boldsymbol {\alpha}_0=(\alpha, 2,\ldots,2)$ and
$$b_{\delta,T}=a_T+ a_T^{-1}\ln \Bigl( 2^{1- m/2}(\Gamma(m/2))^{-1}\mathcal{H}_{D,\alpha}a_T^{2/\alpha+m-2}\Bigr).$$ }
}

\section{Further results and proofs}\label{sec3}

We present first four lemmas followed then by the proofs of \netheo{T1} for $m\ge2$ since the claim for $m=1$,  
\ccL{the stationary} Gaussian processes follows immediately from \cite{HashorvaT2014Different}.
In what follows,  we shall keep the notation as in Section \ref{sec1}, and denote further by $\overline{\Phi}$ and $\varphi$  the survival distribution function and probability density function of a standard normal variable, respectively.  We write $C $ for a positive constant whose values may change from line to line. All the limits are taken as $T$ and $u$ tend to infinity in this coordinated way (unless otherwise stated)
\BQNY
u^{2}=2\ln T+ (2/\alpha+m-2)\ln \ln T+O(1).
\EQNY

Note that, in view of \cite{Pit96},  for any closed non-empty set $ E \subset [0, T ]$
and $\mathcal S_{m-1}$ the unit sphere in $\R^m$ (with respect to $L_2$-norm)
$$\sup_{t\in  E}\chi_{m}(t)=\sup_{(t,\mathbf{v})\in  E\times \mathcal \mathcal S_{m-1}}Y(t,\mathbf{v}),$$
where the Gaussian field $\{Y(t,\mathbf{v}), (t,\mathbf v)\in [0,T]\times \mathcal S_{m-1}\}$ is given by
\BQNY
Y(t,\mathbf{v})=X_{1}(t)v_{1}+X_{2}(t)v_{2}+\cdots+X_{m}(t)v_{m}, \quad (t,\mathbf v)\in [0,T]\times \mathcal S_{m-1}.
\EQNY
Note in passing that the covariance function of $Y(t, \mathbf v)$, denoted by $r(t, \mathbf v, s,\mathbf{w})$, is as follows
\BQN\label{rA}
r(t,\mathbf{v},s,\mathbf{w})=r(t-s)A(\mathbf{v},\mathbf{w}), \quad A(\mathbf{v},\mathbf{w})=1-\frac{\vnorm{\mathbf{v}-\mathbf{w}}^{2}}{2}, \quad\mathbf{v},\mathbf{w}\in \mathcal S_{m-1}.
\EQN
Therefore, crucial in the following is to construct as in \cite{PiterbargS2004Chi} the grids $\mathfrak{R}^\alpha_{b}, b>0$ over the cylinder $[0,T]\times \mathcal S_{m-1}$ (see \eqref{def: RbGrid} for details) and to deal with the random field $Y(t, \mathbf v)$ in terms of $\xi_T(t, \mathbf v)$ defined below in \eqref{def.xi}. \\
Let $\vartheta(x)=\sup_{x\leq |t|\leq T}r(t)$ for any $x>0$. In view of \eqref{corrr}, we choose some small $\ve\in(0, 2^{-1/\alpha})$ such that for all $|t|\leq\varepsilon<2^{-1/\alpha}$
\begin{eqnarray}
\label{eqT20}
\frac{1}{2}|t|^{\alpha}\leq 1-r(t)\leq 2|t|^{\alpha}.
\end{eqnarray}
It follows further from \eqref{cond.Global} that $\vartheta(\varepsilon)<1$ holds for all sufficiently large $T$ (see p.\,86 in \cite{leadbetter1983extremes}). Therefore, we choose some constants $c$ and $a$  such that
\begin{eqnarray}
\label{eqTB}
0<c<a<\frac{1-\vartheta(\varepsilon)}{1+\vartheta(\varepsilon)}<1.
\end{eqnarray}
Next, we  introduce a Gaussian field $\xi_T(t,\mathbf v), (t,\mathbf v)\in [0,T]\times \mathcal S_{m-1}$ via $Y(t,\mathbf v)$ and condition \eqref{cond.Global}, which is crucial in our proof, see the technical \nelem{L3.3}.  Following \cite{Pit2004},  divide $[0,T]$ into intervals with
length $T^{a}$ alternating with shorter intervals with length
$T^{c}$ and write
\BQN\label{def:Interval I,E}
I_i:=
[(i-1)(T^{a}+T^{c}),(i-1)(T^{a}+T^{c})+T^{a}], \quad E_i:=[(i-1)(T^{a}+T^{c}),i(T^{a}+T^{c})),
\EQN
for $1\le i\le n, n=\lfloor T/(T^a+T^c))\rfloor$. Here $\lfloor x\rfloor$ stands for the integer part of $x$. 
\ccL{We will see from \nelem{L3.1} below  that, the asymptotic joint distribution of $(M_m(T), M_m(\delta, T)$ is determined totally by that of the maxima over the closed set $\mathcal I=\cup_{i=1}^n I_i$.} \\
Further, let $Y_{i}(t,\mathbf{v}), (t, \mathbf v)\in [0,T]\times \mathcal S_{m-1}, i\le  n$ be independent copies of  
$\{Y(t,\mathbf{v}), (t, \mathbf v)\in [0,T]\times \mathcal S_{m-1}\}$ 
 \ccL{and} $Z_i, 1\le i\le m$ be standard Gaussian random variables so that
the components of the $(n+m)$-dimension random vector
$$(Y_{1}(t,\mathbf{v}),\ldots,Y_{n}(t,\mathbf{v}),Z_{1},\ldots, Z_{m} )$$ 
are mutually independent. We define, with $\rho(T)=r/\ln T$ and Gaussian random field  $Z(\mathbf{v})=Z_{1}v_{1}+Z_{2}v_{2}+\cdots+Z_{m}v_{m},\ \mathbf{v}\in \mathcal S_{m-1}$,
\BQN\label{def.xi}\xi_{T}(t,\mathbf{v})=\sqrt{1-\rho(T)}Y_{i}(t,\mathbf{v})+\sqrt{\rho(T)}Z(\mathbf{v}), \quad (t,\mathbf{v})\in E_i\times \mathcal S_{m-1}, \ 1\le i\le n,
\EQN
which covariance function 
$\varrho(t, \mathbf{v},s,\mathbf{w})$ 
is given by 
$$\varrho(t, \mathbf{v},s,\mathbf{w})=r^{*}(t,s)A(\mathbf{v},\mathbf{w}),$$
where
\BQN\label{def:r*}
  r^{*}(t,s)=\left\{
 \begin{array}{ll}
  {r(t-s)+(1-r(t-s))\rho(T)},    &(t, s)\in E_i \times E_i;\\
  {\rho(T)},    &(t,s)\in E_i \times E_{j}, i\neq j.
 \end{array}
  \right.
\EQN
\BL\label{L3.2} For the grid
$\mathfrak{R}(\delta)$  is a sparse grid or \ccL{a} Pickands grid, there exists a grid $\mathfrak{R}_{b}^\alpha=\widetilde{\mathfrak{R}}_{b} \times \mathfrak{R}_{b}$ on the cylinder $[0,T]\times \mathcal S_{m-1}$ such that for any
$B>0$, we have for all $x,y\in[-B,B]$
\begin{eqnarray*}
&&\Abs{\pk{a_{T}\big(\max_{t\in \mathcal I }\chi_{m}(t)-b_{T}\big)\leq x,
a_{T}\big(\max_{t\in\mathfrak{R}(\delta)\cap\mathcal I }\chi_{m}(t)- b_{\delta,T}\big)\leq y}\\
&&\ \ \ \ \ -\pk{a_T\big(\max_{(t,\mathbf{v})\in \mathfrak{R}_{b}^{\alpha}\cap(\mathcal I \times \mathcal S_{m-1})}Y(t,\mathbf{v})-b_{T}\big)\leq x,
     a_T\big(\max_{(t,\mathbf{v})\in(\mathfrak{R}(\delta)\times\mathfrak{R}_{b})\cap(\mathcal I \times \mathcal S_{m-1})}Y(t,\mathbf{v})- b_{\delta,T}\big)\leq y }}\rightarrow0
\end{eqnarray*}
as $T\rightarrow\infty$ and $b\downarrow0$\ccL{, subsequently}. 
\EL

For the proof of \nelem{L3.2}, one can follow similar arguments as for Lemma 3 in \cite{PiterbargS2004Chi} and thus we omit here. Since the grid $\mathfrak R_{b}^\alpha$ is crucial for our proofs, we provide the details on its construction. \\
For any given $\ve>0$ we partition the sphere $\mathcal S _{m-1} $ onto $N(\varepsilon)$ parts $A_{1},\ldots,A_{N(\varepsilon)}$ in the
following way. With a polar-coordinate transformation, any point $\mathbf x$ on the sphere $\mathcal S _{m-1}$ is given in terms of angle $\boldsymbol \varphi=(\varphi_{1},\dots, \varphi_{m-1})\in [0,\pi)^{m-2}\times [0,2\pi)$ and
divide \ccL{all} the interval\ccL{s} $[0,\pi]$ into intervals of length $\varepsilon$ (or less for the last
interval), do the same for the interval $[0,2\pi]$. This partition of the parallelepiped $[0,\pi]^{m-2}\times [0,2\pi]$ generates the partition
$A_{j}, 1\le j\le  N(\varepsilon)$ of the sphere. For a fixed $u$, choose in \ccL{every} 
$A_{j}$ an inner point $B_j$ and consider the tangent plane 
to the cylinder $[0,T]\times \mathcal S_{m-1} $ at the chosen point.
Introduce in the tangent plane rectangular coordinates, with origin at the
tangent point; the first coordinate is assigned to the direction $t$.
In the so-constructed spac\ccL{e} 
$\R^{m}$, consider the grid of points
$$\mathfrak{R}_{b,u,\varepsilon}^{j,\alpha,P}:=\left(bl_{1}u^{-\frac{2}{\alpha}}, bl_{2}u^{-1},\ldots,bl_{m} u^{-1}\right),\ \ j=1,2,\ldots,N(\varepsilon)$$
and
$$\mathfrak{R}_{b,u,\varepsilon}^{j,P}:=\left( bl_{2}u^{-1},\ldots,bl_{m} u^{-1}\right),\ \
\mathfrak{\widetilde{R}}_{b,u,\varepsilon}^{j,P}:=\left( bl_{1}u^{-\frac{2}{\alpha}}\right),\,\ j=1,2,\ldots,N(\varepsilon),$$
where $(l_{1},l_{2},\ldots,l _{m} )\in \mathbb Z^{m}$.
Suppose that $\varepsilon$ is so small that the orthogonal projection\ccL{s} of all $[0,T]\times A_{j}$
onto the corresponding tangent plane  are 
one-to-one. \ccL{Hence} 
the distance between any two points in $[0,T]\times \mathcal S_{m-1}$ ha\ccL{s} 
 the same order
with that of their orthogonal projections on the  tangent planes. Denote by $A_{j}^{P}$ the
projection of $A_{j}$ at the tangent plane, and by $\mathfrak{R}_{b,u,\varepsilon}^{j,\alpha}$, $\mathfrak{R}_{b,u,\varepsilon}^{j}$  and $\mathfrak{\widetilde{R}}_{b,u,\varepsilon}^{j}$,
the prototype of $\mathfrak{R}_{b,u,\varepsilon}^{j,\alpha,P}$, $\mathfrak{R}_{b,u,\varepsilon}^{j,P}$
and $\mathfrak{\widetilde{R}}_{b,u,\varepsilon}^{j,P}$, respectively,  under this
projection. The grids
\BQN\label{def: RbGrid}\mathfrak{R}_{b}^{\alpha}=\mathfrak{R}_{b,u,\varepsilon}^{\alpha}=\bigcup_{j=1}^{N(\varepsilon)}\mathfrak{R}_{b,u,\varepsilon}^{j,\alpha},\ \
\mathfrak{R}_{b}=\mathfrak{R}_{b,u,\varepsilon}=\bigcup_{j=1}^{N(\varepsilon)}\mathfrak{R}_{b,u,\varepsilon}^{j},\ \
\mathfrak{\widetilde{R}}_{b}=\mathfrak{\widetilde{R}}_{b,u,\varepsilon}=\bigcup_{j=1}^{N(\varepsilon)}\mathfrak{\widetilde{R}}_{b,u,\varepsilon}^{j}
\EQN
with an appropriate choice of their parameters, satisfy the assertion of \nelem{L3.2}.  \\
Next, we will introduce \ccL{three} technical lemmas which proofs will be relegated in the Appendix. We will see that Lemmas \ref{L3.3} and \ref{L3.4} are crucial for the proof of \netheo{T1}. 
\BL\label{L3.3} Let the grid
$\mathfrak{R}(\delta)$ be  a  sparse grid or  Pickands grid, and $\mathfrak R_b^\alpha$ as in \nelem{L3.2}. For any
$B>0$ we have for all $x,y\in[-B,B]$,
\begin{eqnarray*}
\lefteqn{\Delta_{T,b}:=\bigg|\pk{a_T\big(\max_{(t,\mathbf{v})\in \mathfrak{R}_{b}^{\alpha}\cap(\mathcal I \times \mathcal S_{m-1})}Y(t,\mathbf{v})-b_{T}\big)\leq x,
     a_T\big(\max_{(t,\mathbf{v})\in(\mathfrak{R}(\delta)\times\mathfrak{R}_{b})\cap(\mathcal I \times \mathcal S_{m-1})}Y(t,\mathbf{v})- b_{\delta,T}\big)\leq y}}\notag\\
&&-\pk{a_T\big(\max_{(t,\mathbf{v})\in \mathfrak{R}_{b}^{\alpha}\cap(\mathcal I \times \mathcal S_{m-1})}\xi_{T}(t,\mathbf{v})-b_{T}\big)\leq x,
     a_T\big(\max_{(t,\mathbf{v})\in(\mathfrak{R}(\delta)\times\mathfrak{R}_{b})\cap(\mathcal I \times \mathcal S_{m-1})}\xi_{T}(t,\mathbf{v})- b_{\delta,T}\big)\leq y }\bigg|\rightarrow 0
\end{eqnarray*}
uniformly for $b>0$, as $T\rightarrow\infty$. 
\EL

\COM{\textbf{Lemma 3.4}. {\sl Suppose that the grid
$\mathfrak{R}(\delta)$ is a  sparse grid or  Pickands grid. For any
$B>0$ we have for all $x,y\in[-B,B]$
\begin{eqnarray*}
&&\mathbb P\left\{a_T\big(\max_{(t,\mathbf{v})\in \mathfrak{R}_{b}^{\alpha}\cap(\mathcal I \times \mathcal S_{m-1})}\xi_{T}(t,\mathbf{v})-b_{T}\big)\leq x,
     a_T\big(\max_{(t,\mathbf{v})\in(\mathfrak{R}(\delta)\times\mathfrak{R}_{b})\cap(\mathcal I \times \mathcal S_{m-1})}\xi_{T}(t,\mathbf{v})- b_{\mathbf{T}}^{*}\big)\leq y\right\}\\
&&\ \ \ \ \ \ =\mathbb{E}\prod_{i=1}^{n}
\mathbb P\left\{\max_{(t,\mathbf{v})\in \mathfrak{R}_{b}^{\alpha}\cap(I_i\times \mathcal S_{m-1})}Y_{i}(t,\mathbf{v})\leq v_{T},
   \max_{(t,\mathbf{v})\in(\mathfrak{R}(\delta)\times\mathfrak{R}_{b})\cap(I_i\times \mathcal S_{m-1})}Y_{i}(t,\mathbf{v})\leq v_{T}^{*}\right\},
\end{eqnarray*}
 where
  \begin{gather}
  \begin{aligned} \label{equ413}
&v_{T}:=\frac{b_{T}+x/a_{T}-\sqrt{\rho(T)}z}{(1-\rho(T))^{1/2}}
=\frac{x+r-\sqrt{2r}\vnorm{\mathbf z}}{a_{T}}+b_{T}+o(a_{T}^{-1})\\
v_{T}^{*}:=\frac{b_{\delta,T}+y/a_{T}-\sqrt{\rho(T)}z}{(1-\rho(T))^{1/2}}
=\frac{y+r-\sqrt{2r}\chi_m }{a_{T}}+b_{\delta,T}+o(a_{T}^{-1})
\end{aligned}
\end{gather}
with $b_{\delta,T}$ given by \eqref{def:normalization}.
}
\textbf{Proof:} }
{In the following, we denote (recall $b_{\delta,T}$ in \eqref{def:normalization})
  \begin{gather}
  \begin{aligned} \label{equ413}
&v_{T}:=\frac{b_{T}+x/a_{T}-\sqrt{\rho(T)}\vnorm{\mathbf z}}{(1-\rho(T))^{1/2}}
=\ccL{b_T+}\frac{x+r-\sqrt{2r}\vnorm{\mathbf z}}{a_{T}}+o(a_{T}^{-1})\\
&v_{T}^{*}:=\frac{b_{\delta,T}+y/a_{T}-\sqrt{\rho(T)}\vnorm{\mathbf z}}{(1-\rho(T))^{1/2}}
=\ccL{b_{\delta,T}+}\frac{y+r-\sqrt{2r}\vnorm{\mathbf z}}{a_{T}}+o(a_{T}^{-1}).
\end{aligned}
\end{gather}
}
\BL\label{L3.4}
Under the conditions of \netheo{T1}, 
we have, with $\H_{D,\alpha}, \H_{D, \boldsymbol\alpha_0}^{x,y}, \mathfrak{R}_b$ and $v_{T}, v_{T}^{*}$ given by \eqref{def: H.D}, \eqref{def: RbGrid} and \eqref{equ413}, respectively, 
\BQN\label{LA1.2}
&&\mathbb P\left\{\max_{(t,\mathbf{v})\in \mathfrak{R}_{b}^{\alpha}\cap ([0,T^{a}]\times \mathcal S_{m-1})}Y(t,\mathbf{v})> v_{T}\right\}=T^{a-1}e^{-x-r+\sqrt{2r}\vnorm {\mathbf z}}(1+o(1))\notag \\
&&\mathbb P\left\{\max_{(t,\mathbf{v})\in (\mathfrak{R}(\delta)\times\mathfrak{R}_{b})\cap( [0,T^{a}]\times \mathcal S_{m-1})}Y(t,\mathbf{v})> v_{T}^{*}\right\}
=T^{a-1}e^{-y-r+\sqrt{2r}\vnorm {\mathbf z}}(1+o(1))
\EQN
hold for sufficiently large $T$ and sufficiently small $b>0$. \ccL{And} 
\BQN\label{SparsePickands}
&&\mathbb P\left\{\max_{(t,\mathbf{v})\in \mathfrak{R}_{b}^{\alpha}\cap ([0,T^{a}]\times \mathcal S_{m-1})}Y(t,\mathbf{v})>v_{T},
\max_{(t,\mathbf{v})\in (\mathfrak{R}(\delta)\times\mathfrak{R}_{b})\cap( [0,T^{a}]\times \mathcal S_{m-1})}Y(t,\mathbf{v})> v_{T}^{*}\right\}\notag\\
&&=\left\{\begin{array}{ll}
o(T^{a-1}),& \mathfrak R(\delta) {\rm \ is\ a\ sparse\ grid}; 
 \\
T^{a-1}\pi^{\frac{m-1}2}\mathcal{H}_{D,\boldsymbol {\alpha}_0}^{\ln
 \mathcal{H}_{\alpha}+x, \ln \mathcal{H}_{D,\alpha} +y}e^{-r+\sqrt{2r}\vnorm{\mathbf z}}(1+o(1)), & \mathfrak R(\delta) {\rm \ is\  a\  Pickands\  grid} 
\end{array}
  \right.
\EQN
hold for sufficiently large $T$ and sufficiently small $b>0$.
\EL

\BL\label{L3.1} Suppose that the grid
$\mathfrak{R}(\delta)$ is a sparse grid or \ccL{a} Pickands grid. For any
$B>0$, we have for all $x,y\in[-B,B]$, \Tan{as $T\rightarrow\infty$}
\begin{eqnarray}
\label{eq:L3.1}
&&\Abs{\pk{a_{T}\big(M_{m}(T)-b_{T}\big)\leq x,
           a_{T}\big( M_m(\delta, T)- b_{\delta,T}\big)\leq y}\notag\\
&&\ \ \ \ \ -\pk{a_{T}\big(\max_{t\in \mathcal I }\chi_{m}(t)-b_{T}\big)\leq x,
a_{T}\big(\max_{t\in\mathfrak{R}(\delta)\cap\mathcal I }\chi_{m}(t)- b_{\delta,T}\big)\leq y}}\rightarrow0.
\end{eqnarray}
\EL
\prooflem{L3.1}
The proof is similar to that of Lemma 6 in \cite{Pit2004}. Clearly, the left-hand side of \eqref{eq:L3.1} is bounded from above  by
\begin{eqnarray}
\label{eq400}
\pk{\max_{t\in [0,T]\backslash\mathcal I }\chi_{m}(t)>b_{T}+x/a_{T}}
+\pk{\max_{t\in\mathfrak{R}(\delta)\cap[0,T]\backslash\mathcal I }\chi_{m}(t)>b_{\delta,T}+y/a_{T}}=: J_{T,1}+ J_{T,2}.
\end{eqnarray}
It follows from \ccL{\eqref{def:chi0} and} \eqref{def:normalization.ab} and the construction of $\mathcal I$ that \ccL{(recall that \eqref{def:chi0} holds also for $T=T(u)\to\IF$ with suitable speed, see Theorem 7.2 of \cite{Pit96}),} with $mes(\cdot)$ the Lebesgue measure 
\begin{eqnarray*}
J_{T,1}\le C  mes([0,T]\backslash\mathcal I )(b_{T}+x/a_{T})^{2/\alpha+m-1}\overline{\Phi}(b_{T}+x/a_{T}) \le C  \frac{mes([0,T]\backslash\mathcal I )}{T} \le  C\frac{nT^{c}}{T}\rightarrow0
\end{eqnarray*}
as $T\rightarrow\infty$.
\Tan{Similarly, \ccL{using \eqref{LA1.2} in \nelem{L3.4}  with $v_T^*$ and $\expon{-y-r+\sqrt{2r}\vnorm{\mathbf z}}$ replaced by $u_T^*= b_{\delta,T}+y/a_T$ and $e^{-y}$, respectively,}
} 
we have $\limit{T} J_{T,2}=0$, hence the proof is complete. 
\hfill $\Box$

\textbf{Proof of Theorem 2.1}. First, by \eqref{def.xi} 
we have
\begin{eqnarray}
\label{eq412}
&&\mathbb P\left\{a_T\big(\max_{(t,\mathbf{v})\in \mathfrak{R}_{b}^{\alpha}\cap(\mathcal I \times \mathcal S_{m-1})}\xi_{T}(t,\mathbf{v})-b_{T}\big)\leq x,
     a_T\big(\max_{(t,\mathbf{v})\in(\mathfrak{R}(\delta)\times\mathfrak{R}_{b})\cap(\mathcal I \times \mathcal S_{m-1})}\xi_{T}(t,\mathbf{v})- b_{\delta,T}\big)\leq y\right\}\nonumber\\
&&= \frac{1}{(2\pi)^{m/2}}\int_{\R^{m}}e^{-\frac{1}{2}\vnorm{\mathbf z}^{2}}\mathbb P\left\{a_{T}\big(\max_{(t,\mathbf{v})\in \mathfrak{R}_{b}^{\alpha}\cap(\mathcal I \times \mathcal S_{m-1})}\xi_{T}(t,\mathbf{v})-b_{T}\big)\leq x,\right.\nonumber\\
  &&\ \ \ \left.   a_T\big(\max_{(t,\mathbf{v})\in(\mathfrak{R}(\delta)\times\mathfrak{R}_{b})\cap(\mathcal I \times \mathcal S_{m-1})}\xi_{T}(t,\mathbf{v})- b_{\delta,T}\big)\leq y|Z_{1}=z_{1}, \ccL{\ldots},  Z_{m}=z_{m}\right\}d z_{1}\cdots d z_m \nonumber\\
&&=\int_{\vnorm{\mathbf z}\ge0}\left(\mathbb P\left\{\max_{(t,\mathbf{v})\in \mathfrak{R}_{b}^{\alpha}\cap([0,T^a]\times \mathcal S_{m-1})}Y (t,\mathbf{v})\leq \frac{b_{T}+x/a_{T}-\sqrt{\rho(T)}\vnorm{\mathbf z}}{(1-\rho(T))^{1/2}},\right.\right.\nonumber\\
  &&\ \ \ \left.\left.   \max_{(t,\mathbf{v})\in(\mathfrak{R}(\delta)\times\mathfrak{R}_{b})\cap([0,T^a] \times \mathcal S_{m-1})}Y(t,\mathbf{v})\leq \frac{b_{\delta,T}+y/a_{T}-\sqrt{\rho(T)}\vnorm{\mathbf z}}{(1-\rho(T))^{1/2}}\right\}\right)^n d\pk{\chi_{m}  \le \vnorm{\mathbf z}}. 
\end{eqnarray}
Denote 
\BQN\label{def: Pn}
P_{n,b}(x,y):=\mathbb P\left\{\max_{(t,\mathbf{v})\in\mathfrak{R}_{b}^{\alpha}\cap( [0,T^{a}]\times \mathcal S_{m-1})}Y(t,\mathbf{v})\leq v_{T},
   \max_{(t,\mathbf{v})\in(\mathfrak{R}(\delta)\times\mathfrak{R}_{b})\cap([0,T^a] \times \mathcal S_{m-1})}Y(t,\mathbf{v})\leq v_{T}^{*}\right\}.
\EQN
Next, we deal with the grid \ccL{$\mathfrak{R}(\mathfrak{\delta})$} being a sparse, Pickands and dense grid in turn. \\
\underline{(a) For the sparse grid \ccL{$\mathfrak{R}(\mathfrak{\delta})$}}.
Using Lemmas \ref{L3.2}\ccL{--}\ref{L3.1} and \eqref{eq412}, the first claim in \netheo{T1} will follow if we show that
\begin{eqnarray*}
&& \Abs{(P_{n,b}(x,y))^n-\exp\left(-\big(e^{-x-r+\sqrt{2r}\vnorm{\mathbf z}}+e^{-y-r+\sqrt{2r}\vnorm{\mathbf z}}\big)\right)}\rightarrow 0.
\end{eqnarray*}
Since $\lim_{T\to\IF}P_{n,b}(x,y)=1$ uniformly for all $x,y\in\R$ and thus
$$(P_{n,b}(x,y))^n = \expon{n \ln P_{n,b}(x,y)} = \expon{-n(1- P_{n,b}(x,y)) (1+o(1))}.$$
\COM{\begin{eqnarray*}
&&\prod_{i=1}^{n} \mathbb P\left\{\max_{(t,\mathbf{v})\in \mathfrak{R}_{b}^{\alpha}\cap(I_i\times \mathcal S_{m-1})}Y_{i}(t,\mathbf{v})\leq v_{T},
  \max_{(t,\mathbf{v})\in(\mathfrak{R}(\delta)\times\mathfrak{R}_{b})\cap(I_i\times \mathcal S_{m-1})}Y_{i}(t,\mathbf{v})\leq v_{T}^{*}\right\}\\
&&=\left(\mathbb P\left\{\max_{(t,\mathbf{v})\in\mathfrak{R}_{b}^{\alpha}\cap( [0,T^{a}]\times \mathcal S_{m-1})}Y(t,\mathbf{v})\leq v_{T},
   \max_{(t,\mathbf{v})\in\widetilde{\mathfrak D}}Y(t,\mathbf{v})\leq v_{T}^{*}\right\}\right)^{n}\\
&&=\exp\left(n\ln\left(\mathbb P\left\{\max_{(t,\mathbf{v})\in \mathfrak{R}_{b}^{\alpha}\cap ([0,T^{a}]\times \mathcal S_{m-1})}Y(t,\mathbf{v})\leq v_{T},
   \max_{(t,\mathbf{v})\in\widetilde{\mathfrak D}}Y(t,\mathbf{v})\leq v_{T}^{*}\right\}\right)\right)\\
&&=\exp\left(-n\left(1-\mathbb P\left\{\max_{(t,\mathbf{v})\in \mathfrak{R}_{b}^{\alpha}\cap ([0,T^{a}]\times \mathcal S_{m-1})}Y(t,\mathbf{v})\leq v_{T},
   \max_{(t,\mathbf{v})\in\widetilde{\mathfrak D}}Y(t,\mathbf{v})\leq v_{T}^{*}\right\}\right)+R_{n}\right).
\end{eqnarray*}
Since $P_{n,b}(x,y)\to1$ 
we get that the remainder
$R_{n}$ can be estimated as
$R_{n}=o(n(1-P_{n,b}(x,y)))$. }
Finally, using \nelem{L3.4} for sparse grids, we get that
\begin{eqnarray*}
n\left(1-P_{n,b}(x,y)
\right)= nT^{a-1}\left(e^{-x-r+\sqrt{2r}\vnorm {\mathbf z}}+e^{-y-r+\sqrt{2r}\vnorm {\mathbf z}}\right) (1+o(1)), 
\end{eqnarray*}
which together with the fact that $n=T/(T^a+T^c), 0<c<a<1$ and the dominated convergence theorem completes the
proof for sparse grid. \\
\underline{(b) For the Pickands grid \ccL{$\mathfrak{R}(\mathfrak{\delta})$} with $\delta(T)=D(2\ln T)^{-1/\alpha}$}. 
Similarly as for the sparse grid,  it suffices to show that
\begin{eqnarray*}
&&n\left(1-P_{n,b}(x,y) \right)=\ccL{\left(e^{-x}+e^{-y}
       -\pi^{(m-1)/2}\H_{a,\boldsymbol {\alpha}_0}^{\ln \H_{a,\alpha}+x,\ln \H_{\alpha}+y}\right)e^{-r+\sqrt{2r}\vnorm{\mathbf z}}} (1+o(1))
\end{eqnarray*}
with $P_{n,b}(x,y)$ defined in \eqref{def: Pn}. 
\COM{\begin{eqnarray*}
&&n\left(1-\mathbb P\left\{\max_{(t,\mathbf{v})\in \mathfrak{R}_{b}^{\alpha}\cap ([0,T^{a}]\times \mathcal S_{m-1})}Y(t,\mathbf{v})\leq v_{T},
   \max_{(t,\mathbf{v})\in\widetilde{\mathfrak D}}Y(t,\mathbf{v})\le v_{T}^{*}\right\}\right)\\
&&=n\mathbb P\left\{\max_{(t,\mathbf{v})\in \mathfrak{R}_{b}^{\alpha}\cap ([0,T^{a}]\times \mathcal S_{m-1})}Y(t,\mathbf{v})> v_{T}\right\}
   +n\mathbb P\left\{\max_{(t,\mathbf{v})\in\widetilde{\mathfrak D}}Y(t,\mathbf{v})> v_{T}^{*}\right\}\\
&&\ \ -n\mathbb P\left\{\max_{(t,\mathbf{v})\in \mathfrak{R}_{b}^{\alpha}\cap ([0,T^{a}]\times \mathcal S_{m-1})}Y(t,\mathbf{v})> v_{T},
   \max_{(t,\mathbf{v})\in\widetilde{\mathfrak D}}Y(t,\mathbf{v})> v_{T}^{*}\right\}\\
&&= nT^{a-1}\left(e^{-x-r+\sqrt{2r}\chi_m }+e^{-y-r+\sqrt{2r}\chi_m }- \pi^{(m-1)/2}\mathcal{H}_{D,\boldsymbol {\alpha}_0}^{\ln
 \mathcal{H}_{\alpha}+x, \ln \mathcal{H}_{D,\alpha} +y}e^{-r+\sqrt{2r}\chi_m }\right)(1+o(1)), 
\end{eqnarray*}
}
This is verified by 
\nelem{L3.4} for Pickands grids.\\ 
\underline{(c) For the dense grid $\mathfrak{R}(\delta)$}.
In view of Lemma 3 of \cite{PiterbargS2004Chi}, we have
\begin{eqnarray*}
&&\bigg|\mathbb P\left\{a_{T}\big(\max_{(t,\mathbf{v})\in [0,T]\times \mathcal S_{m-1}}Y(t,\mathbf{v})-b_{T}\big)\leq x, a_{T}\big(\max_{(t,\mathbf{v})\in \mathfrak{R}(\delta)\cap[0,T]\times \mathcal S_{m-1}}Y(t,\mathbf{v})- b_{T}\big)\leq y\right\}\\
&&\ \ \ -\mathbb P\left\{a_{T}\big(\max_{(t,\mathbf{v})\in [0,T]\times \mathcal S_{m-1}}Y(t,\mathbf{v})-b_{T}\big)\leq x, a_{T}\big(\max_{(t,\mathbf{v})\in [0,T]\times \mathcal S_{m-1}}Y(t,\mathbf{v})- b_{T}\big)\leq y\right\}\bigg|\\
&&\leq\bigg|\mathbb P\left\{a_{T}\big(\max_{(t,\mathbf{v})\in \mathfrak{R}(\delta)\cap[0,T]\times \mathcal S_{m-1}}Y(t,\mathbf{v})- b_{T}\big)\leq y\right\}
-\mathbb P\left\{a_{T}\big(\max_{(t,\mathbf{v})\in [0,T]\times \mathcal S_{m-1}}Y(t,\mathbf{v})- b_{T}\big)\leq y\right\}\bigg| \rightarrow 0.
\end{eqnarray*}
Further, by \eqref{eq TB}
\BQNY
&&\mathbb P\left\{a_{T}\big(\max_{(t,\mathbf{v})\in [0,T]\times \mathcal S_{m-1}}Y(t,\mathbf{v})-b_{T}\big)\leq x, a_{T}\big(\max_{(t,\mathbf{v})\in [0,T]\times \mathcal S_{m-1}}Y(t,\mathbf{v})- b_{T}\big)\leq y\right\} \\
&&=\mathbb P\left\{a_{T}\big(\max_{(t,\mathbf{v})\in [0,T]\times \mathcal S_{m-1}}Y(t,\mathbf{v})-b_{T}\big)\leq \min(x,y) \right\}\\
&& \rightarrow \mathbb{E}\exp\left(-e^{-\min(x,y)-r+\sqrt{2r}\chi_m }\right),
\EQNY
the required claim (c) follows. Consequently, \netheo{T1} is proved. \QED
\COM{\textbf{Proof of Theorem 2.3}. First, note that
\begin{eqnarray*}
&&\mathbb P\left\{a_{T}\big(M_{m}(T)-b_{T}\big)\leq x, a_{T}\big( M_m(\delta, T)- b_{T}\big)\leq y\right\}\\
&&=\mathbb P\left\{a_{T}\big(\max_{(t,\mathbf{v})\in [0,T]\times \mathcal S_{m-1}}Y(t,\mathbf{v})-b_{T}\big)\leq x, a_{T}\big(\max_{(t,\mathbf{v})\in \mathfrak{R}(\delta)\cap[0,T]\times \mathcal S_{m-1}}Y(t,\mathbf{v})- b_{T}\big)\leq y\right\}.
\end{eqnarray*}
In view of Lemma 3 of \cite{PiterbargS2004Chi} we have
\begin{eqnarray*}
&&\bigg|\mathbb P\left\{a_{T}\big(\max_{(t,\mathbf{v})\in [0,T]\times \mathcal S_{m-1}}Y(t,\mathbf{v})-b_{T}\big)\leq x, a_{T}\big(\max_{(t,\mathbf{v})\in \mathfrak{R}(\delta)\cap[0,T]\times \mathcal S_{m-1}}Y(t,\mathbf{v})- b_{T}\big)\leq y\right\}\\
&&\ \ \ -\mathbb P\left\{a_{T}\big(\max_{(t,\mathbf{v})\in [0,T]\times \mathcal S_{m-1}}Y(t,\mathbf{v})-b_{T}\big)\leq x, a_{T}\big(\max_{(t,\mathbf{v})\in [0,T]\times \mathcal S_{m-1}}Y(t,\mathbf{v})- b_{T}\big)\leq y\right\}\bigg|\\
&&\leq\bigg|\mathbb P\left\{a_{T}\big(\max_{(t,\mathbf{v})\in \mathfrak{R}(\delta)\cap[0,T]\times \mathcal S_{m-1}}Y(t,\mathbf{v})- b_{T}\big)\leq y\right\}
-\mathbb P\left\{a_{T}\big(\max_{(t,\mathbf{v})\in [0,T]\times \mathcal S_{m-1}}Y(t,\mathbf{v})- b_{T}\big)\leq y\right\}\bigg| \rightarrow 0,
\end{eqnarray*}
as $T\rightarrow\infty.$
Next, applying  Theorem 2.1, we get
\begin{eqnarray*}
&&\mathbb P\left\{a_{T}\big(\max_{(t,\mathbf{v})\in [0,T]\times \mathcal S_{m-1}}Y(t,\mathbf{v})-b_{T}\big)\leq x, a_{T}\big(\max_{(t,\mathbf{v})\in [0,T]\times \mathcal S_{m-1}}Y(t,\mathbf{v})- b_{T}\big)\leq y\right\}\\
&&=P\{a_{T}\big(\max_{(t,\mathbf{v})\in [0,T]\times \mathcal S_{m-1}}Y(t,\mathbf{v})-b_{T}\big)\leq \min(x, y)\}\\
&&\rightarrow \mathbb{E}\exp\left(-e^{-\min(x,y)-r+\sqrt{2r}\chi_m }\right),
\end{eqnarray*}
as $T\rightarrow\infty$,
hence the proof is complete. \hfill$\Box$
}

\section{Appendix} \label{sec4}

In this section, we give the proofs of Lemmas \ref{L3.3} and \ref{L3.4}, respectively. 
Before we proceed the proof, let us recall some basic quantities which will be repeatedly used below.  For simplicity of notation, we write $u_{T}=
b_{T}+x/a_T, u_{T}^{*}=
b_{\delta,T}+y/a_{T}$ with $a_T, b_T, b_{\delta, T}$ given by \eqref{def:normalization.ab} and \eqref{def:normalization}. Thus
\begin{eqnarray}
\label{eq403}
u^{2}_{T}=2\ln T+ 
(2/\alpha+m-2)\ln \ln T+O(1)\ccL,
\end{eqnarray}
which implies that
\BQN\label{uT}
T^{-1} =C u^{2/\alpha+m-2}\expon{-\frac{u^2}2}(1+o(1)).
\EQN
Further, denote by $\varpi(t,s)=\max\{|r(t-s)|, |r^{*}(t,s)|\}$ with $r^*$ given in \eqref{def:r*}, and define
$$
\theta(t_0)=\sup_{0\leq t,s\leq T,|t-s|>t_0} \varpi(t,s), \quad t_0>0. 
$$
Since $\theta(t_0)\ge \vartheta(t_0):=\sup_{t_0\le \abs t\le  T} r(t)$, the constants $c$ and $a$ given in \eqref{eqTB} hold also for $\theta(\cdot)$, i.e.,
\begin{eqnarray*}
0<c<a<\frac{1-\theta(\varepsilon)}{1+\theta(\varepsilon)}<1.
\end{eqnarray*}
Note that, from the construction of $\mathfrak{R}_{b}$ given by \eqref{def: RbGrid}, the number of points in $\mathfrak{R}_{b}\cap \mathcal S_{m-1}$ does not
exceed $\cL{Cb^{-(m-1)}}u^{m-1}(1+o(1))=C  b^{-(m-1)}(2\ln T)^{(m-1)/2}(1+o(1))$.
\\
\prooflem{L3.3} It follows by Berman's inequality  (
see e.g., \cite{Pit96}) that, with  $\mathcal D_i= \mathfrak{R}_{b}^{\alpha}\cap(I_i \times \mathcal S_{m-1}), \widetilde {\mathcal D}_i= (\mathfrak{R}(\delta)\times\mathfrak{R}_{b})\cap(I_i \times \mathcal S_{m-1}),  i\le n$ 
\begin{eqnarray}
\label{decomp.Delta}
\Delta_{T,b} &\leq& \sum_{(t,\mathbf{v})\in \mathfrak D_i, (s,\mathbf{w})\in \mathfrak D_j,\atop (t,\mathbf{v})\neq (s,\mathbf{w}), 1\leq i,j\leq n }|\Upsilon_{r,\varrho}|
\int_{0}^{1}\frac{1}{\sqrt{1-r^{(h)}(t,\mathbf{v},s,\mathbf{w})}}\exp\left(-\frac{u_{T}^{2}}{1+r^{(h)}(t,\mathbf{v},s,\mathbf{w})}\right)dh\nonumber\\
&&\quad+\sum_{(t,\mathbf{v})\in \widetilde{\mathfrak D}_i, (s,\mathbf{w})\in \widetilde{\mathfrak D}_j,\atop (t,\mathbf{v})\neq (s,\mathbf{w}), 1\leq i,j\leq n }|\Upsilon_{r,\varrho}|
\int_{0}^{1}\frac{1}{\sqrt{1-r^{(h)}(t,\mathbf{v},s,\mathbf{w})}}\exp\left(-\frac{(u_{T}^{*})^{2}}{1+r^{(h)}(t,\mathbf{v},s,\mathbf{w})}\right)dh\nonumber\\
&&\quad+\sum_{(t,\mathbf{v})\in \mathfrak D_i, (s,\mathbf{w})\in \widetilde{\mathfrak D}_j, \atop (t,\mathbf{v})\neq (s,\mathbf{w}), 1\leq i,j\leq n }|\Upsilon_{r,\varrho}|
\int_{0}^{1}\frac{1}{\sqrt{1-r^{(h)}(t,\mathbf{v},s,\mathbf{w})}}\exp\left(-\frac{u_{T}^{2}+(u_{T}^{*})^{2}}{2(1+r^{(h)}(t,\mathbf{v},s,\mathbf{w}))}\right)dh\nonumber\\
&=:& \Delta_{T,b}^{(1)} +\Delta_{T,b}^{(2)} +\Delta_{T,b}^{(3)}\ccL{,}
\end{eqnarray}
where $\Upsilon_{r,\varrho}=r(t,\mathbf{v},s,\mathbf{w})-\varrho(t,\mathbf{v},s,\mathbf{w})$ and
$r^{(h)}(t,\mathbf{v},s,\mathbf{w})=hr(t,\mathbf{v},s,\mathbf{w})+(1-h)\varrho(t,\mathbf{v},s,\mathbf{w})$. Next, we shall show that $\Delta_{T,b}^{(i)} =o(1), i=1,2,3$ for sufficiently large $T $ and small $b>0$, respectively. \\
We shall present first the proof for $\Delta_{T,b}^{(1)} =o(1)$. To this end, we consider below the sum with $(t,\mathbf v), (s, \mathbf w)$  in the same $\mathfrak D_i, 1\leq i \leq n$, denoted by $\Delta_{T,b}^{(1,1)}$,
 and split further the sum into two parts as follows 
\begin{eqnarray}
\label{eq402}
\Delta_{T,b}^{(1,1)}=\sum_{(t,\mathbf{v}), (s,\mathbf{w})\in \mathcal D_i,\atop 1\le i\le n,|t-s|\leq \varepsilon }
+\sum_{(t,\mathbf{v}), (s,\mathbf{w})\in \mathcal D_i,\atop 1\le i\le n,|t-s|> \varepsilon }=:J_{T,1}+J_{T,2}
\end{eqnarray}
for some small $\ve>0$ such that \eqref{eqT20} and \eqref{eqTB} hold.
For $J_{{T},1}$, note that in this case, it follows from \eqref{def:r*} that 
$\abs{\Upsilon_{r,\varrho}}=\rho(T)(1-r(t-s)) A(\mathbf v, \mathbf w)$,  and
by (\ref{corrr}) that, we can choose small  $\varepsilon>0$ such that
\BQNY
r^{(h)}(t,\mathbf{v},s,\mathbf{w})=(r (t-s)+(1-h)(1-r(t-s))\rho(T))A(\mathbf v, \mathbf w)
= r(t-s)A(\mathbf v, \mathbf w) (1+o(1)) \EQNY
holds for sufficiently large $T$ and $|t-s|\leq\varepsilon$.
Consequently, we have (recall that $A(\mathbf{v},\mathbf{w})\leq 1$)
\begin{eqnarray*}
\label{eq404}
J_{T,1}&\leq&C \sum_{{(t,\mathbf{v}), (s,\mathbf{w})\in \mathcal D_i,\atop 1\le i\le n,|t-s|\leq \varepsilon }}{\rho(T)\sqrt{1-r(t-s)}}\exp\left(-\frac{u^{2}_{T}}{1+r(t-s)|A(\mathbf{v},\mathbf{w})|}\right)\nonumber\\
&\leq &\ccL{C}Tb^{-1}u_{T}^{2/\alpha}\rho(T)\exp\left(-\frac{u^{2}_{T}}{2}\right)\sum_{\mathbf{v}, \mathbf{w}\in \mathfrak{R}_{b}\cap  \mathcal S_{m-1},\atop t\in \mathcal{\widetilde{R}}_{b}\cap[0,T], |t|\leq \varepsilon }\sqrt{1-r(t)}
\exp\left(-\frac{(1-r(t)|A(\mathbf{v},\mathbf{w})|)u^{2}_{T}}{2(1+r(t)|A(\mathbf{v},\mathbf{w})|)}\right)\nonumber\\
&= &C Tb^{-1}u_{T}^{2/\alpha}\rho(T)\exp\left(-\frac{u^{2}_{T}}{2}\right)\sum_{\mathbf{v}, \mathbf{w}\in \mathfrak{R}_{b}\cap  \mathcal S_{m-1},\atop t\in \mathcal{\widetilde{R}}_{b}\cap[0,T],|t|\leq \varepsilon }\sqrt{1-r(t)}
\exp\left(-\frac{(1-r(t))u^{2}_{T}}{2(1+r(t))}\right)\nonumber\\
&&\ \ \times \exp\left(-\frac{r(t)(1-|A(\mathbf{v},\mathbf{w})|)u^{2}_{T}}{2(1+r(t))(1+r(t)|A(\mathbf{v},\mathbf{w})|)}\right)\nonumber\\
&\leq & C  Tb^{-m}u_{T}^{2/\alpha}u_{T}^{m-1}\rho(T)\exp\left(-\frac{u^{2}_{T}}{2}\right)\sum_{\mathbf{v}\in \mathfrak{R}_{b}\cap  \mathcal S_{m-1},\atop t\in \mathcal{\widetilde{R}}_{b}\cap[0,T], |t|\leq \varepsilon }\sqrt{1-r(t)}
\exp\left(-\frac{(1-r(t))u^{2}_{T}}{2(1+r(t))}\right)\nonumber\\
&&\ \ \times \exp\left(-\frac{r(t)(1-|A(\mathbf{v},\mathbf{w}_{0})|)u^{2}_{T}}{2(1+r(t))(1+r(t)|A(\mathbf{v},\mathbf{w}_{0})|)}\right),
\end{eqnarray*}
where $\mathbf{w}_{0}$ is any fixed point on $ \mathfrak{R}_{b} \cap \mathcal S_{m-1} $. Since
\begin{eqnarray*}
\lefteqn{\sum_{\mathbf{v}\in \mathfrak{R}_{b}\cap  \mathcal S_{m-1},\atop t\in \mathcal{\widetilde{R}}_{b}\cap[0,T],|t|\leq \varepsilon }
\exp\left(-\frac{r(t)(1-|A(\mathbf{v},\mathbf{w}_{0})|)u^{2}_{T}}{2(1+r(t))(1+r(t)|A(\mathbf{v},\mathbf{w}_{0})|)}\right)}\\
&&\leq \sum_{\mathbf{v}\in \mathfrak{R}_{b}\cap \mathcal S_{m-1}}\exp\left(-C u_{T}^{2}\vnorm{\mathbf{v}-\mathbf{w}_{0}}^{2}\right)\leq C,
\end{eqnarray*}
it follows further by (\ref{eqT20}), \eqref{uT} and $\rho(T)=r/\ln T=O(u_{T}^{-2})$ that
\begin{eqnarray*}
J_{T,1}&\leq & C  Tb^{-m}u_{T}^{2/\alpha}u_{T}^{m-1}\rho(T)\exp\left(-\frac{u^{2}_{T}}{2}\right)\sum_{t\in \mathcal{\widetilde{R}}_{b}\cap[0,T],|t|\leq \varepsilon }\sqrt{1-r(t)}
\exp\left(-\frac{(1-r(t))u^{2}_{T}}{2(1+r(t))}\right)\nonumber\\
&\leq & C  b^{-m}u_{T}^{-1}\sum_{ t\in \mathcal{\widetilde{R}}_{b}\cap[0,T],|t|\leq \varepsilon }\sqrt{2}|t|^{\alpha/2}\exp\left(-\frac{|t|^{\alpha}u^{2}_{T}}{8}\right)\nonumber\\
&\leq & C  b^{-m}u_{T}^{-1}\sum_{k=1}^{\infty}\expon{-\frac{1}{4}(kb)^{\alpha}}\nonumber\\
&\leq & C  b^{-m}u_{T}^{-1},
\end{eqnarray*}
which implies that $J_{T,1}=o(1)$ uniformly for $b>0$ as $T\rightarrow\infty$.\\
Using the fact that $u_{T}=a_T (1+o(1))$, we obtain
\begin{eqnarray}
\label{eq405}
J_{T,2}&\leq&C \sum_{{(t,\mathbf{v}), (s,\mathbf{w})\in \mathcal D_i,\atop 1\le i\le n,|t-s|>\varepsilon }}
\exp\left(-\frac{u^{2}_{T}}{1+\cL{\abs{r(t-s)}}}\right)\nonumber\\
&\leq & C T^{1+a}b^{-2m}u_{T}^{4/\alpha}u_{T}^{2m-2}
     \exp\left(-\frac{u^{2}_{T}}{1+\theta(\varepsilon)}\right)\nonumber\\
&\leq & C  T^{1+a}b^{-2m}u_{T}^{4/\alpha}u_{T}^{2m-2} T^{-\frac{2}{1+\theta(\varepsilon)}}\nonumber\\
&\leq & C   T^{a-\frac{1-\theta(\varepsilon)}{1+\theta(\varepsilon)}}b^{-2m}(\ln T)^{2/\alpha+m-1}.
\end{eqnarray}
Thus, $J_{T,2}=o(1)$  uniformly for $b>0$ as $T\to\IF$ since $a<(1-\theta(\varepsilon))/({1+\theta(\varepsilon)})$.

Next, we consider the sum $\Delta_{T,b}^{(1)}$ with $(t,\mathbf v), (s, \mathbf w)$  in $\mathfrak D_i, \mathfrak D_j$ with $1\leq i\neq j \leq n$, denoted by $\Delta_{T,b}^{(1,0)}$.
Note that in this case, $|t-s|\geq T^{c}$ and $\varrho(s,\mathbf{v},t,\mathbf{w})=\rho(T)A(\mathbf{v},\mathbf{w})$. 
 Choose $\beta$ such that
$0<c<a<\beta <(1-\theta(\varepsilon))/(1+\theta(\varepsilon))$ and
split the sum$\Delta_{T,b}^{(1,0)}$ into two parts as follows
\begin{eqnarray}
\label{eq402.2}
\Delta_{T,b}^{(1,0)}= \sum_{(t,\mathbf{v})\in \mathcal D_i, (s,\mathbf{w})\in  \mathcal D_j, \atop1\leq i\neq j\leq n, |t-s|\leq T^{\beta} }
+\sum_{(t,\mathbf{v})\in \mathcal D_i, (s,\mathbf{w})\in  \mathcal D_j, \atop1\leq i\neq j\leq n, |t-s|> T^{\beta} }=:S_{T,1}+S_{T,2}.
\end{eqnarray}
For $S_{T,1}$, with the similar derivation as for
(\ref{eq405}), we have
\begin{eqnarray}
\label{eq407}
S_{T,1}&\leq&C \sum_{(t,\mathbf{v})\in \mathcal D_i, (s,\mathbf{w})\in  \mathcal D_j, \atop1\leq i\neq j\leq n, |t-s|\leq T^{\beta} }
\exp\left(-\frac{u^{2}_{T}}{1+r(t-s)}\right)\nonumber\\
&\leq & C T^{1+\beta}b^{-2m}u_{T}^{4/\alpha}u_{T}^{2m-2}
     \exp\left(-\frac{u^{2}_{T}}{1+\theta(\varepsilon)}\right)\nonumber\\
&\leq & C  T^{1+\beta}b^{-2m}u_{T}^{4/\alpha}u_{T}^{2m-2} T^{-\frac{2}{1+\theta(\varepsilon)}}\nonumber\\
&\leq & C   T^{\beta-\frac{1-\theta(\varepsilon)}{1+\theta(\varepsilon)}}b^{-2m}(\ln T)^{2/\alpha+m-1}
\end{eqnarray}
implying that $S_{T,1}=o(1)$ uniformly for $b>0$, since $\beta<(1-\theta(\varepsilon))/({1+\theta(\varepsilon)})$.
\\
For $S_{T,2}$, we need \ccL{some} more precise estimation.
By condition (\ref{cond.Global}), there exists some constant $ K>0$
such that
$$\theta(t)\ln t\leqslant K$$
for $t, T$ sufficiently large. 
Thus $\theta(t)\leq K/\ln T^{\beta}, \  t\ge T^\beta$ holds for $T$ large enough.
Now using \eqref{eq403}, we 
obtain
\begin{eqnarray}
\label{Tan2}
T^{2}u_{T}^{4/\alpha+2m-2}(\ln T)^{-1}\exp\left(-\frac{u_{T}^{2}}{1+\theta(T^{\beta})}\right)&\leq& T^{2}u_{T}^{4/\alpha+2m-2}(\ln T)^{-1}\exp\left(-\frac{u_{T}^{2}}{1+K/\ln T^{\beta}}\right)\nonumber\\
&\leq& C   \left(T^{2}(\ln T)^{2/\alpha+m-2}\right)^{1-\frac{1}{1+K/\ln T^{\beta}}} \leq C.
\end{eqnarray}
Therefore, by similar arguments as for Lemma 6.4.1 of
\cite{leadbetter1983extremes} we have 
\begin{eqnarray*}
S_{T,2}&\leq&C \sum_{(t,\mathbf{v})\in \mathcal D_i, (s,\mathbf{w})\in  \mathcal D_j, \atop1\leq i\neq j\leq n, |t-s|> T^{\beta} }
      |r(t-s)-\rho(T)|\exp\left(-\frac{u^{2}_{T}}{1+\theta(T^{\beta})}\right)\nonumber\\
&\leq&C Tb^{-\cL{(2m-\ccL{2})}}u_{T}^{2/\alpha}u_{T}^{2m-2}\exp\left(-\frac{u^{2}_{T}}{1+\theta(T^{\beta})}\right)
      \sum_{t\in \mathfrak{\widetilde{R}}_{b}\cap[0,T], t>T^{\beta}}|r(t)-\rho(T)|\nonumber\\
&=& C T^{2}(\ln T)^{-1}u_{T}^{4/\alpha}u_{T}^{2m-2}\exp\left(-\frac{u^{2}_{T}}{1+\theta(T^{\beta})}\right)\cdot b^{-\cL{(2m-\ccL{2})}}\frac{\ln T}{Tu_{T}^{2/\alpha}}
     \sum_{t\in \mathfrak{\widetilde{R}}_{b}\cap[0,T], t>T^{\beta}}|r(t)-\rho(T)|\nonumber\\
&\leq& C b^{-\cL{(2m-\ccL{2})}}\frac{\ln T}{Tu_{T}^{2/\alpha}}
     \sum_{t\in \mathfrak{\widetilde{R}}_{b}\cap[0,T], t>T^{\beta}}|r(t)-\rho(T)|\nonumber\\
&\leq& C b^{-\cL{(2m-\ccL{2})}}\frac{1}{\beta Tu_{T}^{2/\alpha}}
     \sum_{t\in \mathfrak{\widetilde{R}}_{b}\cap[0,T], t>T^{\beta}}|r(t)\ln t-r|+ C b^{-\cL{(2m-\ccL{2})}}\frac{r}{Tu_{T}^{2/\alpha}}
     \sum_{t\in \mathfrak{\widetilde{R}}_{b}\cap[0,T], t>T^{\beta}}\bigg|1-\frac{\ln T}{\ln t}\bigg|,
\end{eqnarray*}
where, by \eqref{cond.Global}, the first term is  $o(1)$  uniformly for $b>0$, and the second term is also $o(1)$  uniformly for $b>0$ following an integral estimate below (see also the proof of Lemma 6.4.1 in \cite{leadbetter1983extremes})
\begin{eqnarray*}
C b^{-\cL{(2m-1)}}\frac{r}{Tu_{T}^{2/\alpha}}
     \sum_{t\in \mathfrak{\widetilde{R}}_{b}\cap[0,T], t>T^{\beta}}\Abs{1-\frac{\ln T}{\ln t}}
&\leq& C b^{-\cL{(2m-1)}}\frac{r}{Tu_{T}^{2/\alpha}}\frac{1}{\ln T^{\beta}}
     \sum_{t\in \mathfrak{\widetilde{R}}_{b}\cap[0,T], t>T^{\beta}}\abs{\ln t-\ln T}\\
&=& C  b^{-\cL{(2m-1)}}\frac{r}{\ln T^{\beta}}\int_0^1 |\ln x|dx.
\end{eqnarray*}
Consequently, combining the assertions for $J_{T,i}, S_{T,i}, i=1,2$ in \eqref{eq402}, \eqref{eq402.2}, we have $\Delta_{T,b}^{(1)}= o(1)$. \\
The proof of $\Delta_{T,b}^{(2)}=o(1)$ is similar  as that for $\Delta_{T,b}^{(1)}=o(1)$ with minor modifications by replacing $u_T, \mathfrak D_i, i\le n$ by $u_T^*, \widetilde{\mathfrak D}_i, i\le n$, we omit thus the details.\\
It remains to prove $\Delta_{T,b}^{(3)}=o(1)$.
Recall that $\mathfrak{R}(\delta)$ can be a sparse grid or \ccL{a} Pickands grid. We only show below the proof for $\mathfrak{R}(\delta)$ a sparse grid by  following the main arguments as for $\Delta_{T,b}^{(1)}$. The Pickands grid case can be shown similarly for the sparse grid and thus we omit \ccL{it} here. \\
Consider first the sum $\Delta_{T,b}^{(3)}$ with $t,s$ in the
same $I_i, i\le n$, which is further split  into two parts as
\begin{eqnarray}
\label{eqB31}
\Delta_{T,b}^{(3,1)}:=\sum_{(t,\mathbf{v})\in \mathcal D_i, (s,\mathbf{w})\in \widetilde{ \mathcal D}_i, \atop  1\leq i \leq n, |t-s|\leq \varepsilon }
+\sum_{(t,\mathbf{v})\in \mathcal D_i, (s,\mathbf{w})\in \widetilde{ \mathcal D}_i, \atop  1\leq i \leq n, |t-s|>\varepsilon }=:\widetilde J_{T,1}+\widetilde J_{T,2}.
\end{eqnarray}
Note that 
\begin{eqnarray}
\label{eqB32}
\widetilde u_T^2:=\frac{1}{2}(u_{T}^{2}+(u_{T}^{*})^{2})=2\ln T+\ln a_{T}^{2/\alpha+m-2}+\ln (\delta^{-1}a_{T}^{m-2})+O(1)
\end{eqnarray}
and the grid $\mathfrak R(\delta)$ is a sparse grid, i.e., $\limit{T}\delta (2\ln T)^{1/2}=\IF$, the remaining proof of $\widetilde J_{T,1}+ \widetilde J_{T,2}=o(1)$ is similar to that for $J_{T,1}$ and thus we omit it here.
\COM{
Consequently, we have
\begin{eqnarray*}
\label{eq404}
\widetilde J_{T,1}&\leq&C \sum_{(t,\mathbf{v})\in \mathcal D_i, (s,\mathbf{w})\in \widetilde{ \mathcal D}_i, \atop  1\leq i \leq n, |t-s|\leq \varepsilon }
|r(t-s)-r^{*}(t,s)|\frac{1}{\sqrt{1-r(t-s)}}\exp\left(-\frac{w^{2}_{T}}{1+r(t-s)|A(\mathbf{v},\mathbf{w})|}\right)\nonumber\\
&\leq &\mathcal{ C}T\delta^{-1}\rho(T)\exp\left(-\frac{w^{2}_{T}}{2}\right)\sum_{\mathbf{v}, \mathbf{w}\in \mathfrak{R}_{b}\cap  \mathcal S_{m-1},\atop t\in \mathfrak{\widetilde{R}}_{b}\cap[0,T], |t|\leq \varepsilon }\sqrt{1-r(t)}
\exp\left(-\frac{(1-r(t)|A(\mathbf{v},\mathbf{w})|)w^{2}_{T}}{2(1+r(t)|A(\mathbf{v},\mathbf{w})|)}\right)\nonumber\\
&\leq & C  T\delta^{-1}\rho(T)\exp\left(-\frac{w^{2}_{T}}{2}\right)\sum_{\mathbf{v}, \mathbf{w}\in \mathfrak{R}_{b}\cap  \mathcal S_{m-1},\atop t\in \mathfrak{\widetilde{R}}_{b}\cap[0,T],|t|\leq \varepsilon }\sqrt{1-r(t)}
\exp\left(-\frac{(1-r(t))w^{2}_{T}}{2(1+r(t))}\right)\nonumber\\
&&\ \ \times \exp\left(-\frac{r(t)(1-|A(\mathbf{v},\mathbf{w})|)w^{2}_{T}}{2(1+r(t))(1+r(t)|A(\mathbf{v},\mathbf{w})|)}\right)\nonumber\\
&\leq & C  T\delta^{-1}b^{-(m-1)}u_{T}^{m-1}\rho(T)\exp\left(-\frac{w^{2}_{T}}{2}\right)\sum_{\mathbf{v}\in \mathfrak{R}_{b}\cap  \mathcal S_{m-1},\atop t\in \mathfrak{\widetilde{R}}_{b}\cap[0,T], |t|\leq \varepsilon }\sqrt{1-r(t)}
\exp\left(-\frac{(1-r(t))w^{2}_{T}}{2(1+r(t))}\right)\nonumber\\
&&\ \ \times \exp\left(-\frac{r(t)(1-|A(\mathbf{v},\mathbf{w}_{0})|)w^{2}_{T}}{2(1+r(t))(1+r(t)|A(\mathbf{v},\mathbf{w}_{0})|)}\right),
\end{eqnarray*}
where $\mathbf{w}_{0}$ is a fixed point on $ \mathfrak{R}_{b} \cap \mathcal S_{m-1} $. Since
\begin{eqnarray*}
&&\sum_{\mathbf{v}\in \mathfrak{R}_{b}\cap  \mathcal S_{m-1},\atop t\in \mathfrak{\widetilde{R}}_{b}\cap[0,T],|t|\leq \varepsilon }
\exp\left(-\frac{r(t)(1-|A(\mathbf{v},\mathbf{w}_{0})|)w^{2}_{T}}{2(1+r(t))(1+r(t)|A(\mathbf{v},\mathbf{w}_{0})|)}\right)\\
&&\leq \sum_{\mathbf{v}\in \mathfrak{R}_{b}\cap  \mathcal S_{m-1}}\exp\left(-C \widetilde u_T^2\vnorm{\mathbf{v}-\mathbf{w}_{0}}^{2}\right)\leq C
\end{eqnarray*}
then by (\ref{eqT20}) and (\ref{eqB32})
\begin{eqnarray*}
\widetilde J_{T,1}&\leq & C  T\delta^{-1}b^{-(m-1)}u_{T}^{m-1}\rho(T)\exp\left(-\frac{w^{2}_{T}}{2}\right)\sum_{ t\in \mathfrak{\widetilde{R}}_{b}\cap[0,T],|t|\leq \varepsilon }\sqrt{1-r(t)}
\exp\left(-\frac{(1-r(t))w^{2}_{T}}{2(1+r(t))}\right)\nonumber\\
&\leq & C b^{-(m-1)} [(\ln T)^{1/\alpha+1} \delta]^{-1/2}\sum_{ t\in \mathfrak{\widetilde{R}}_{b}\cap[0,T],|t|\leq \varepsilon }\sqrt{2}|t|^{\alpha/2}\exp\left(-\frac{|t|^{\alpha}w^{2}_{T}}{8}\right)\nonumber\\
&\leq & C b^{-(m-1)} [(\ln T)^{1/\alpha+1} \delta]^{-1/2}\sum_{k=1}^{\infty}\expon{-\frac{1}{4}(kb)^{\alpha}}\nonumber\\
&\leq & C b^{-(m-1)} [(\ln T)^{1/\alpha+1} \delta]^{-1/2},
\end{eqnarray*}
which shows $\widetilde J_{T,1}\rightarrow0$ uniformly for $b>0$ as $T\rightarrow\infty$, since $(\ln T)^{1/\alpha}\delta\rightarrow\infty$ for sparse grid.\\
Using the fact that $w_{T}\thicksim (2\ln
T)^{1/2}$, we obtain
\begin{eqnarray*}
\widetilde J_{T,2}&\leq&C \sum_{(t,\mathbf{v})\in \mathcal D_i, (s,\mathbf{w})\in \widetilde{ \mathcal D}_i, \atop  1\leq i \leq n, |t-s|> \varepsilon }
\exp\left(-\frac{w^{2}_{T}}{1+r(t-s)}\right)\nonumber\\
&\leq & C T^{1+a}u_{T}^{2/\alpha}\delta^{-1}b^{-\cL{2(m-1)}}u_{T}^{2m-2}
     \exp\left(-\frac{w^{2}_{T}}{1+\theta(\varepsilon)}\right)\nonumber\\
&\leq & C  T^{1+a}u_{T}^{2/\alpha}\delta^{-1}b^{-\cL{2(m-1)}}u_{T}^{2m-2} T^{-\frac{2}{1+\theta(\varepsilon)}}\nonumber\\
&\leq & C T^{a-\frac{1-\theta(\varepsilon)}{1+\theta(\varepsilon)}}b^{-\cL{2(m-1)}}(\ln T)^{2/\alpha+m-1}.
\end{eqnarray*}
Thus, $\widetilde J_{T,2}\rightarrow0$ uniformly for $b>0$ as $T\rightarrow\infty$, since $a<\frac{1-\theta(\varepsilon)}{1+\theta(\varepsilon)}$.}

\COM{Second, as for $\Delta_{T,b}^{(1)}$, we consider the case that $t,s$ in the
different intervals $  I_i$ and $I_j$ for $i\neq j$.
We also have, in this case,
$|t-s|\geq T^{c}$ and $\varrho(s,\mathbf{v},t,\mathbf{w})=\rho(T)A(\mathbf{v},\mathbf{w})$ for  $(s,\mathbf{v})\in   I_i\times \mathcal S_{m-1}$
and $(t,\mathbf{w})\in I_j\times \mathcal S_{m-1}$, $i\neq j$. Choose $\beta$ such that
$0<c<a<\beta<\frac{1-\theta(\varepsilon)}{1+\theta(\varepsilon)}$ and
split the sum (\ref{eqB11}) into two parts as
\begin{eqnarray*}
\sum_{(t,\mathbf{v})\in \mathcal D_i, (s,\mathbf{w})\in \widetilde{ \mathcal D}_j, \atop  1\leq i\neq j\leq n, |t-s|\leq T^{\beta} }
+\sum_{(t,\mathbf{v})\in \mathcal D_i, (s,\mathbf{w})\in \widetilde{ \mathcal D}_j, \atop  1\leq i\neq j\leq n, |t-s|>T^{\beta} }=:M_{T,1}+M_{T,2}.
\end{eqnarray*}
For $M_{T,1}$, with the similar derivation as for
$\widetilde J_{T,2}$, we have
\begin{eqnarray*}
\label{eq407}
M_{T,1}&\leq&C \sum_{(t,\mathbf{v})\in \mathcal D_i, (s,\mathbf{w})\in \widetilde{ \mathcal D}_j, \atop  1\leq i\neq j\leq n, |t-s|\leq T^{\beta} }
\exp\left(-\frac{w^{2}_{T}}{1+r(t-s)}\right)\nonumber\\
&\leq & C T^{1+\beta}u_{T}^{2/\alpha}\delta^{-1}b^{-\cL{2(m-1)}}u_{T}^{2m-2}
     \exp\left(-\frac{w^{2}_{T}}{1+\theta(\varepsilon)}\right)\nonumber\\
&\leq & C  T^{1+\beta}u_{T}^{2/\alpha}\delta^{-1}b^{-\cL{2(m-1)}}u_{T}^{2m-2}T^{-\frac{2}{1+\theta(\varepsilon)}}\nonumber\\
&\leq & C T^{\beta-\frac{1-\theta(\varepsilon)}{1+\theta(\varepsilon)}}b^{-\cL{2(m-1)}}(\ln T)^{2/\alpha+m-1}.
\end{eqnarray*}
Consequently, since $\beta<\frac{1-\theta(\varepsilon)}{1+\theta(\varepsilon)}$, we have $M_{T,1}\rightarrow0$ uniformly for $b>0$ as $T\rightarrow\infty$.\\
For bound $M_{T,2}$, we use the estimation of the term $S_{T,2}$  with some changes.
As for \eqref{uT}, 
using (\ref{eqB32}), we can  also show that
\begin{eqnarray*}
&&T^{2}u_{T}^{2/\alpha+2m-2}\delta^{-1}(\ln T)^{-1}\exp\left(-\frac{\widetilde u_T^2}{1+\theta(T^{\beta})}\right)=O(1).
\end{eqnarray*}
Therefore, we have
\begin{eqnarray*}
\label{eq509}
M_{T,2}&\leq&C \sum_{(t,\mathbf{v})\in \mathcal D_i, (s,\mathbf{w})\in \widetilde{ \mathcal D}_j, \atop  1\leq i\neq j\leq n, |t-s|> T^{\beta} }
      |r(t-s)-\rho(T)|\exp\left(-\frac{w^{2}_{T}}{1+\theta(T^{\beta})}\right)\nonumber\\
&\leq&C T\delta^{-1}b^{-2(m-1)}u_{T}^{2m-2}\exp\left(-\frac{w^{2}_{T}}{1+\theta(T^{\beta})}\right)
      \sum_{t\in \mathfrak{\widetilde{R}}_{b}\cap[0,T], t>T^{\beta}}|r(t)-\rho(T)|\nonumber\\
&=& C T^{2}(\ln T)^{-1}\delta^{-1}u_{T}^{2/\alpha}u_{T}^{2m-2}\exp\left(-\frac{w^{2}_{T}}{1+\theta(T^{\beta})}\right)\cdot b^{-2(m-1)}\frac{\ln T}{Tu_{T}^{2/\alpha}}
     \sum_{t\in \mathfrak{\widetilde{R}}_{b}\cap[0,T], t>T^{\beta}}|r(t)-\rho(T)|\nonumber\\
&\leq& C b^{-\cL{2(m-1)}}\frac{\ln T}{Tu_{T}^{2/\alpha}}
     \sum_{t\in \mathfrak{\widetilde{R}}_{b}\cap[0,T], t>T^{\beta}}|r(t)-\rho(T)|.
\end{eqnarray*}
The remainder proof is the same as that of $S_{T,2}$, so we omit it.}
Next, for the remaining sum $\Delta_{T,b}^{(3)}-\Delta_{T,b}^{(3,1)}$, i.e., the summand with $t,s$ in the different intervals $I_i, I_j,  1\le i\neq j\le n$, one can show that  (recall  (\ref{eqB32}))
\begin{eqnarray*}
&&T^{2}u_{T}^{2/\alpha+2m-2}\delta^{-1}(\ln T)^{-1}\exp\left(-\frac{\widetilde u_T^2}{1+\theta(T^{\beta})}\right)=O(1).
\end{eqnarray*}
The rest proof of $\Delta_{T,b}^{(3)}-\Delta_{T,b}^{(3,1)}=o(1)$ is the same as that for \eqref{eq402.2}. Consequently, we have $\Delta_{T,b}^{(3)}=o(1)$, which together with   \eqref{decomp.Delta} and the proved $\Delta_{T,b}^{(i)}=o(1), i=1,2$, completes  the proof of \nelem{L3.3}. \hfill$\Box$
\\
\prooflem{L3.4} First, noting that
\begin{eqnarray*}
\mathbb P\left\{\max_{(t,\mathbf{v})\in \mathfrak{R}_{b}^{\alpha}\cap ([0,T^{a}]\times \mathcal S_{m-1})}Y(t,\mathbf{v})> v_{T}\right\}
\cL{\leq} \mathbb P\left\{\max_{(t,\mathbf{v})\in [0,T^{a}]\times \mathcal 
S_{m-1}}Y(t,\mathbf{v})> v_{T}\right\} = \pk{M_m(T^a) >v_T},
\end{eqnarray*}
the first assertion follows thus by \eqref{def:chi0} with elementary calculations, \cL{since \eqref{def:chi0} holds also for $T=T(u)\to\IF$ with suitable speed, see Theorem 7.2 of \cite{Pit96}.} \\
Next, we will show the proofs of \eqref{LA1.2} and \eqref{SparsePickands} with $\mathfrak R(\delta)$ a Pickands and sparse grid, respectively. \\
\underline{Proof of \eqref{LA1.2} with $\mathfrak R(\delta)$ a Pickands  or sparse grid}.  The proof for the Pickands grid is similar as that for Corollary 7.3 in \cite{Pit96} with minor modification (replacing $\H_\alpha$ by $\H_{D, \alpha}$), and thus we omit the details here. \\
Now, we consider $\mathfrak R(\delta)$ a sparse grid. For simplicity, we denote in the following $\cL{\widetilde{\mathcal D} =(\mathfrak{R}  (\delta)\times\mathfrak{R}_{b})\cap( [0,T^{a}]\times \mathcal S_{m-1})}$ and
\BQNY
P_{T,b}^{(1)}:=\sum_{t\in \mathfrak{R}(\delta)\cap[0,T^{a}]}\mathbb{P}\left\{\max_{\mathbf{v}\in\mathfrak{R}_{b}\cap  \mathcal S_{m-1}}Y(t,\mathbf{v})> v_{T}^{*}\right\},\quad P_{T,b}^{(2)}:=\sum_{(t,\mathbf{v})\neq
(s,\mathbf{w})\in\widetilde{\mathfrak D}}\mathbb{P}\left\{Y(t,\mathbf{v})> v_{T}^{*},Y(s,\mathbf{w})> v_{T}^{*}\right\}.
\EQNY
By Bonferroni's inequality, we have 
\begin{eqnarray*}
 P_{T,b}^{(1)}-P_{T,b}^{(2)} \leq \mathbb P\left\{\max_{(t,\mathbf{v})\in\widetilde{\mathfrak D}}Y(t,\mathbf{v})> v_{T}^{*}\right\} \leq  P_{T,b}^{(1)},
\end{eqnarray*}
therefore, it suffices to show that
$$P_{T,b}^{(1)}=T^{a-1}e^{-y-r+\sqrt{2r}\vnorm {\mathbf z}}(1+o(1)), \quad P_{T,b}^{(2)} =o(P_{T,b}^{(1)})
$$
hold \cL{for sufficiently large $T$ and small $b>0$}. Clearly,
\begin{eqnarray*}
P_{T,b}^{(1)}&=&(1+o(1))\sum_{t\in\mathfrak{R}(\delta)\cap[0,T^{a}]}\mathbb{P}\left\{\max_{\mathbf{v}\in \mathcal S_{m-1}}Y(t,\mathbf{v})> v_{T}^{*}\right\} \\
&=&(1+o(1))T^{a}\delta^{-1}\pk{\chi_m(t) > v_{T}^{*}}\\
&=&(1+o(1)) T^{a-1}e^{-y-r+\sqrt{2r}\vnorm {\mathbf z}}
\end{eqnarray*}
following by elementary calculations. 
It remains to deal with $P_{T,b}^{(2)}$. Split the term $P_{T,b}^{(2)}$ into two parts as
\begin{eqnarray}
\label{decom:P}
P_{T,b}^{(2)}= \sum_{(t,\mathbf{v})\neq (s,\mathbf{w})\in\widetilde{\mathfrak D}, |t-s|<\epsilon}
+\sum_{(t,\mathbf{v})\neq
(s,\mathbf{w})\in\widetilde{\mathfrak D}, |t-s|\geq \epsilon}
=: P_{T,21}+P_{T,22}.
\end{eqnarray}
Using the well-known results for bivariate Gaussian tail probability (see e.g., p.\,225 in \cite{leadbetter1983extremes}), we have
\begin{eqnarray*}
\label{eq.A12}
P_{T,21}&\le &\sum_{(t,\mathbf{v})\neq
(s,\mathbf{w})\in\widetilde{\mathfrak D}, |t-s|<\epsilon}\left[\overline{\Phi}(v_{T}^{*})
\overline{\Phi}\left(v_{T}^{*}\frac{\sqrt{1-r(t,\mathbf{v},s,\mathbf{w})}}{\sqrt{1+r(t,\mathbf{v},s,\mathbf{w})}}\right)\right].
\end{eqnarray*}
By \eqref{rA} and \eqref{eqT20}, we can choose $\epsilon>0$ small enough such that
$$\frac{1-r(t,\mathbf{v},s,\mathbf{w})}{1+r(t,\mathbf{v},s,\mathbf{w})}\geq \frac{1}{4}|t-s|^{\alpha}
+\frac{1}{8}\vnorm{\mathbf{v}-\mathbf{w}}^{2}$$
and we thus have
\begin{eqnarray*}
\label{eq.A12}
P_{T,21}&\leq&C \sum_{(t,\mathbf{v})\neq
(s,\mathbf{w})\in\widetilde{\mathfrak D}, |t-s|<\epsilon}\left[\overline{\Phi}(v_{T}^{*})
\overline{\Phi}\left(v_{T}^{*}\sqrt{\frac{|t-s|^{\alpha}}4
+\frac{\vnorm{\mathbf{v}-\mathbf{w}}^{2}}8}\right)\right]\\
&\leq& C \overline{\Phi}(v_{T}^{*})\sum_{(t,\mathbf{v})\neq
(s,\mathbf{w})\in\widetilde{\mathfrak D}, |t-s|<\epsilon}
\frac{1}{|t-s|^{\alpha/2}v_{T}^{*}}\exp\left(-\frac{1}{8}|t-s|^{\alpha}(v_{T}^{*})^{2}\right)
\exp\left(-\frac{1}{16}\vnorm{\mathbf{v}-\mathbf{w}}^{2}(v_{T}^{*})^{2}\right)\\
&=& C \overline{\Phi}(v_{T}^{*})\sum_{t, s\in\mathfrak{R}(\delta)\cap [0,T^{a}], 0<|t-s|<\epsilon}
\frac{1}{|t-s|^{\alpha/2}v_{T}^{*}}\exp\left(-\frac{1}{8}|t-s|^{\alpha}(v_{T}^{*})^{2}\right)\\
&&\times \sum_{\mathbf{v}\neq
\mathbf{w}\in\mathfrak{R}_{b}\cap \mathcal S_{m-1}}
\exp\left(-\frac{1}{16}\vnorm{\mathbf{v}-\mathbf{w}}^{2}(v_{T}^{*})^{2}\right)\\
&\leq& C \overline{\Phi}(v_{T}^{*})\sum_{t,
s\in\mathfrak{R}(\delta)\cap [0,T^{a}], 0<|t-s|<\epsilon}
\frac{1}{|t-s|^{\alpha/2}v_{T}^{*}}\exp\left(-\frac{1}{8}|t-s|^{\alpha}(v_{T}^{*})^{2}\right)\\
&&\times b^{-(m-1)}u_{T}^{m-1}\sum_{\mathbf{v}\in\mathfrak{R}_{b}\cap \mathcal S_{m-1}}
\exp\left(-\frac{1}{16}\vnorm{\mathbf{v}-\mathbf{w}_0}^{2}(v_{T}^{*})^{2}\right).
\end{eqnarray*}
where $\mathbf{w}_{0}$ is any fixed point on $\mathfrak{R}_{b}\cap \mathcal S_{m-1}$. Since
\begin{eqnarray*}
&&\sum_{\mathbf{v}\in\mathfrak{R}_{b}\cap \mathcal S_{m-1}}
\exp\left(-\frac{1}{16}\vnorm{\mathbf{v}-\mathbf{w}_{0}}^{2}(v_{T}^{*})^{2}\right)\leq  C,
\end{eqnarray*}
for sufficiently large $T$. 
Using further the definition of $v_{T}^{*}$ we obtain
\begin{eqnarray*}
P_{T,21}&\leq & C T^{a}\delta^{-1}b^{-(m-1)}u_{T}^{m-1}\overline{\Phi}(v_{T}^{*})\sum_{0< k\delta\leq \epsilon}\frac{1}{(k\delta)^{\alpha/2}v_{T}^{*}}\exp\left(-\frac{1}{8}(k\delta)^{\alpha}(v_{T}^{*})^{2}\right)\\
&=&C  T^{a-1}b^{-(m-1)}\sum_{0< k\delta\leq \epsilon}\frac{1}{[k\delta(\ln T)^{1/\alpha}]^{\alpha/2}}\exp\left(-\frac{1}{4}[k\delta(\ln T)^{1/\alpha}]^{\alpha}\right)(1+o(1))\\
&\leq &C  T^{a-1}b^{-(m-1)}\frac{1}{[(\ln T)^{1/\alpha}\delta]^{\alpha/2}}\sum_{0<k\leq \lfloor\epsilon/\delta\rfloor+1}\exp\left(-\frac{1}{4}[k\delta(\ln T)^{1/\alpha}]^{\alpha}\right)(1+o(1))\\
&\leq &C  T^{a-1}b^{-(m-1)}\frac{1}{[(\ln T)^{1/\alpha}\delta]^{\alpha/2}}(1+o(1))\\
&=&T^{a-1}b^{-(m-1)}o(1),
\end{eqnarray*}
where we used additionally  the fact that $\lim_{T\rightarrow\infty}(\ln T)^{1/\alpha}\delta=\infty$, since $\mathfrak{R}(\delta)$ is a sparse grid. Thus,
we have $P_{T,21}=o(T^{a-1})$ uniformly for $b>0$ as $T\rightarrow\infty$.

For the second term $P_{T,22}$ in \eqref{decom:P}, by Normal Comparison Lemma, we have
\begin{eqnarray*}
P_{T,22}
&\leq & \sum_{(t,\mathbf{v})\neq
(s,\mathbf{w})\in\widetilde{\mathfrak D}, |t-s|\geq \epsilon}\bigg[\overline{\Phi}^{2}(v_{T}^{*})
+C \exp\left(-\frac{(v_{T}^{*})^{2}}{1+|r(t,\mathbf{v},s,\mathbf{w})|}\right)\bigg]\\
&\leq &C T^{a}\delta^{-1}b^{-2(m-1)}u_{T}^{2(m-1)} \sum_{\epsilon\leq k\delta\leq T^{a}}\bigg[\overline{\Phi}^{2}(v_{T}^{*})
+C \exp\left(-\frac{(v_{T}^{*})^{2}}{1+|r(k\delta)|}\right)\bigg]\\
&\leq & C T^{2a}\delta^{-2}b^{-2(m-1)}u_{T}^{2(m-1)}\bigg[\overline{\Phi}^{2}(v_{T}^{*})
+C \exp\left(-\frac{(v_{T}^{*})^{2}}{1+\vartheta(\epsilon)}\right)\bigg]\\
&=:&P_{T,221}+P_{T,222}.
\end{eqnarray*}
By \eqref{equ413}, we have 
$v_{T}^{*}= u_{T}(1+o(1))$. Therefore,
\begin{eqnarray*}
P_{T,221}&\leq& C T^{2a}\delta^{-2}b^{-2(m-1)}u_{T}^{2(m-1)}\frac{\varphi^{2}(v_{T}^{*})}{(v_{T}^{*})^{2}}\\
&\leq&  C T^{2a}\delta^{-2}b^{-2(m-1)}u_{T}^{2(m-2)}\exp\left(-(v_{T}^{*})^{2}\right)\\
&\leq& C T^{2a}\delta^{-2}b^{-2(m-1)}u_{T}^{2(m-2)}[T^{-1}\delta u_{T}^{-(m-2)}]^{2}\\
&=&o(T^{a-1})
\end{eqnarray*}
uniformly for $b>0$ as $T\rightarrow\infty$.
Since  $u_{T}=v_{T}^{*}(1+o(1))= (2\ln T)^{1/2}(1+o(1))$
\begin{eqnarray*}
\label{eq.A16}
P_{T,222}&\leq&C T^{2a}\delta^{-2}b^{-2(m-1)}u_{T}^{2(m-1)}\exp\left(-\frac{(v_{T}^{*})^{2}}{1+\vartheta(\epsilon)}\right)\\
&\leq & C T^{2a}\delta^{-2}b^{-2(m-1)}u_{T}^{2(m-1)}T^{-\frac{2}{1+\vartheta(\epsilon)}}\nonumber\\
&\leq & C  T^{a-1} T^{a-\frac{1-\vartheta(\epsilon)}{1+\vartheta(\epsilon)}}\delta^{-2}b^{-2(m-1)}(\ln T)^{m-1}.
\end{eqnarray*}
Both (\ref{eqTB}) and  $(\ln T)^{1/\alpha}\delta =\infty$ imply
$ S_{T,22}=o(T^{a-1})$ uniformly for $b>0$ as $T\rightarrow\infty$. This completes the proof of the second assertion.

\underline{Proof of \eqref{SparsePickands} with $\mathfrak R(\delta)$ a sparse grid}
\COM{ First, repeating the proof of Lemma 3 of \cite{PiterbargS2004Chi}, we can show
\begin{eqnarray*}
&&\bigg|\mathbb P\left\{\max_{(t,\mathbf{v})\in [0,T^{a}]\times \mathcal S_{m-1}}Y(t,\mathbf{v})>v_{T},
\max_{(t,\mathbf{v})\in\widetilde{\mathfrak D}}Y(t,\mathbf{v})> v_{T}^{*}\right\}\\
&&\ \ -\mathbb P\left\{\max_{(t,\mathbf{v})\in \mathfrak{R}_{b}^{\alpha}\cap ([0,T^{a}]\times \mathcal S_{m-1})}Y(t,\mathbf{v})>v_{T},
\max_{(t,\mathbf{v})\in\widetilde{\mathfrak D}}Y(t,\mathbf{v})> v_{T}^{*}\right\}\bigg|=C T^{a-1}g(b)
\end{eqnarray*}
as $T\rightarrow\infty$, where $g(b)\rightarrow0$ and $b\downarrow0$. So in order to prove the third assertion, it suffices to show
\begin{eqnarray}
\label{eq.A11}
\mathbb P\left\{\max_{(t,\mathbf{v})\in \mathfrak{R}_{b}^{\alpha}\cap ([0,T^{a}]\times \mathcal S_{m-1})}Y(t,\mathbf{v})>v_{T},
\max_{(t,\mathbf{v})\in\widetilde{\mathfrak D}}Y(t,\mathbf{v})> v_{T}^{*}\right\}=o(T^{a-1})
\end{eqnarray}
uniformly for $b>0$ as $T\rightarrow\infty$.
Note that the two maxima in (\ref{eq.A11}) may have the same points, but it does not affect the
 probability. So, without loss of generality, we assume that they do not contain the same points.}
 For simplicity, we denote below $\mathfrak D:=\mathfrak{R}_{b}^{\alpha}\cap ([0,T^{a}]\times \mathcal S_{m-1})$.
Obviously, we have
\begin{eqnarray*}
&&\mathbb P\left\{\max_{(t,\mathbf{v})\in \mathfrak D}Y(t,\mathbf{v})>v_{T},
\max_{(t,\mathbf{v})\in\widetilde{\mathfrak D}}Y(t,\mathbf{v})> v_{T}^{*}\right\}\\
&&\cL=\sum_{(t,\mathbf{v})\in\mathfrak D, (s,\mathbf{w})\in\widetilde{\mathfrak D}, |t-s|<\epsilon}+\sum_{(t,\mathbf{v})\in\mathfrak D, (s,\mathbf{w})\in\widetilde{\mathfrak D}, |t-s| \geq \epsilon}=: Q_{T,21}+Q_{T,22}.
\end{eqnarray*}
By the same argument as for the term $P_{T,21}$, we have for $\mathbf{w}_{0}$ fixed on $\mathfrak{R}_{b}\cap \mathcal S_{m-1}$
\begin{eqnarray*}
\label{eq.A12}
Q_{T,21}&\leq&C \sum_{(t,\mathbf{v})\in\mathcal D, (s,\mathbf{w})\in\widetilde{\mathfrak D}, |t-s|<\epsilon}\left[\overline{\Phi}(v_{T})
\overline{\Phi}\left(v_{T}^{*}(\frac{1}{4}|t-s|^{\alpha}
+\frac{1}{8}\vnorm{\mathbf{v}-\mathbf{w}}^{2})^{1/2}\right)\right]\\
&\leq& C \overline{\Phi}(v_{T})\sum_{(t,\mathbf{v})\in\mathcal D,
(s,\mathbf{w})\in\widetilde{\mathfrak D}, |t-s|<\epsilon}
\frac{1}{|t-s|^{\alpha/2}v_{T}^{*}}\exp\left(-\frac{1}{8}|t-s|^{\alpha}(v_{T}^{*})^{2}\right)
\exp\left(-\frac{1}{16}\vnorm{\mathbf{v}-\mathbf{w}}^{2}(v_{T}^{*})^{2}\right)\\
&\leq& C \overline{\Phi}(v_{T})\sum_{t\in\mathfrak{\widetilde{R}}_{b}\cap [0,T^{a}]\atop
s\in\mathfrak{R}(\delta)\cap [0,T^{a}], |t-s|<\epsilon}
\frac{1}{|t-s|^{\alpha/2}v_{T}^{*}}\exp\left(-\frac{1}{8}|t-s|^{\alpha}(v_{T}^{*})^{2}\right)\\
&&\times b^{-(m-1)}u_{T}^{(m-1)}\sum_{\mathbf{v}\in\mathfrak{R}_{b}\cap \mathcal S_{m-1}}
\exp\left(-\frac{1}{16}\vnorm{\mathbf{v}-\mathbf{w}_{0}}^{2}(v_{T}^{*})^{2}\right)\\
&\leq & C T^{a}b^{\ccL1-m}u_{T}^{2/\alpha}u_{T}^{(m-1)}\overline{\Phi}(v_{T})\sum_{0< k\delta\leq \epsilon}\frac{1}{(k\delta)^{\alpha/2}v_{T}^{*}}\exp\left(-\frac{1}{8}(k\delta)^{\alpha}(v_{T}^{*})^{2}\right)\\
&\leq&C  T^{a-1}b^{\ccL1-m}\sum_{0< k\delta\leq \epsilon}\frac{1}{(k\delta)^{\alpha/2}(\ln T)^{1/2}}\exp\left(-\frac{1}{4}(k\delta)^{\alpha}\ln T\right)\\
&\leq &C  T^{a-1}b^{\ccL1-m}\frac{1}{(\ln T)^{1/2}\delta^{\alpha/2}}\sum_{0<k\leq [\epsilon/\delta]+1}
\exp\left(-\frac{1}{4}(k\delta)^{\alpha}\ln T\right)\\
&\leq &C  T^{a-1}b^{\ccL1-m}\frac{1}{[(\ln T)^{1/\alpha}\delta]^{\alpha/2}}\\
&=&T^{a-1}o(1),
\end{eqnarray*}
uniformly for $b>0$, where we used additionally  the fact that $\lim_{T\rightarrow\infty}(\ln T)^{1/\alpha}\delta=\infty$,
since $\mathfrak{R}(\delta)$ is a  sparse grid.

To bound the term $Q_{T,22}$, using again Normal Comparison Lemma, with the same arguments as for the term $P_{T,22}$,
 we have
\begin{eqnarray*}
Q_{T,22}
&\leq & \sum_{(t,\mathbf{v}), (s,\mathbf{w})\in\widetilde{\mathfrak D}, |t-s|\geq \epsilon}\bigg[\overline{\Phi}(v_{T})\overline{\Phi}(v_{T}^{*})
+C \exp\left(-\frac{v_{T}^2+(v_{T}^{*})^{2}}{2(1+|r(t,\mathbf{v},s,\mathbf{w})|)}\right)\bigg]\\
&\leq &C T^{a}\delta^{-1}b^{-2(m-1)}u_{T}^{2(m-1)} 
\sum_{\epsilon\leq k\delta\leq T^{a}}\bigg[\overline{\Phi}(v_{T})\overline{\Phi}(v_{T}^{*})
+C \exp\left(-\frac{ v_{T}^2+(v_{T}^{*})^{2}}{2(1+|r(k\delta)|)}\right)\bigg]\\
&\leq & C T^{2a}\delta^{-2}b^{-\cL{2(m}-1)}u_{T}^{2(m-1)} 
\bigg[\overline{\Phi}(v_{T})\overline{\Phi}(v_{T}^{*})
+C \exp\left(-\frac{ v_{T}^2+(v_{T}^{*})^{2}}{2(1+\vartheta(\epsilon))}\right)\bigg]\\
&=:&Q_{T,221}+Q_{T,222}.
\end{eqnarray*}
By the same arguments as for $P_{T,221}$ and $P_{T,222}$,
we can show that $Q_{T,221}=o(T^{a-1})$ and $Q_{T,222}=o(T^{a-1})$ uniformly for $b>0$ as $T\rightarrow\infty$, respectively. Consequently,  \eqref{SparsePickands} holds for $\mathfrak R(\delta)$ a sparse grid. 

\COM{\BL\label{L3.4}
Under the conditions of Theorem 2.2
we have
$$\mathbb P\left\{\max_{(t,\mathbf{v})\in \mathfrak{R}_{b}^{\alpha}\cap ([0,T^{a}]\times \mathcal S_{m-1})}Y(t,\mathbf{v})> v_{T}\right\}=T^{a-1}e^{-x-r+\sqrt{2r}\vnorm{\mathbf z}}(1+o(1)),$$
as $b\downarrow0$ and $T\rightarrow\infty$;
$$\mathbb P\left\{\max_{(t,\mathbf{v})\in\widetilde{\mathfrak D}}Y(t,\mathbf{v})> v_{T}^{*}\right\}
=T^{a-1}e^{-y-r+\sqrt{2r}\vnorm{\mathbf z}}(1+o(1))$$
as $b\downarrow0$ and $T\rightarrow\infty$, and $\H_{D, \boldsymbol\alpha}^{x,y}\in(0,\IF)$
\begin{eqnarray*}
&&\mathbb P\left\{\max_{(t,\mathbf{v})\in \mathfrak{R}_{b}^{\alpha}\cap ([0,T^{a}]\times \mathcal S_{m-1})}Y(t,\mathbf{v})>v_{T},
\max_{(t,\mathbf{v})\in\widetilde{\mathfrak D}}Y(t,\mathbf{v})> v_{T}^{*}\right\}\\
&&\ \ \ \ \ =T^{a-1}\pi^{(m-1)/2}\mathcal{H}_{D,\boldsymbol {\alpha}_0}^{\ln
 \mathcal{H}_{\alpha}+x, \ln \mathcal{H}_{D,\alpha} +y}e^{-r+\sqrt{2r}\vnorm{\mathbf z}}(1+o(1)),
\end{eqnarray*}
as $b\downarrow0$ and $T\rightarrow\infty$, where $v_{T}$ and $v_{T}^{*}$ are defined in \eqref{equ413}. 
\EL

\textbf{Proof:} The first assertion is just the one in Lemma A1 and we present it here just for citing easily. Recall that $\mathfrak{R}(\delta)=\mathfrak{R}(D(2\ln
T)^{-1/\alpha})$ with $D>0$  is a Pickands grid under the conditions of Theorem 2.2.
Thus, let first $b\downarrow0$ and then the second assertion can be proved by repeating the proof of Corollary 7.3 of \cite{Pit96} by replacing the Pickands constant $\mathcal{H}_{\alpha}$
by $\mathcal{H}_{D,\alpha}$. We omit the details.
}

\underline{Proof of \eqref{SparsePickands} with $\mathfrak{R}(\delta)$ a Pickands grid}.
We shall use below some notation and results from \cite{Pit96}. Let $\mathcal{A}$ be a set in $\R^{m}$ and $\mathbf{d}=(d_1, \ldots, d_m)$ with all $d_i>0,i\le m$; denote
$$\mathbf{d}\mathcal{A}=\bigg(\mathbf{x}=(x_{1},x_{2},\ccL{\ldots,} x_{m}):\left(\frac{x_{1}}{d_{1}},\frac{x_{2}}{d_{2}}, \ccL{\ldots,} \frac{x_{m}}{d_{m}}\right)\in \mathcal{A}\bigg)$$
and with  $\lambda>0$  a constant
$$\mathbf{g}_{u}=(u^{-2/\alpha_{1}},u^{-2/\alpha_{2}}, \ccL{\ldots}, u^{-2/\alpha_{m}}),\quad \mathcal{K}=[0,\lambda]^{m}.$$
Let $Z(\mathbf{t}), \mathbf{t}\in \R^{m}$ be a homogeneous Gaussian random field with correlation function \ccL{$r_{Z}(\mathbf{t})$  such that}, for some $\alpha_{i}\in(0,2], i\le m$
$$r_{Z}(\mathbf{t})=1-\sum_{i=1}^{m}|t_{i}|^{\alpha_{i}}(1+o(1)), \ \ \vnorm{\mathbf t}\rightarrow 0\ \ \mbox{and}\ \ r_{Z}(\mathbf{t})<1,\ \ \forall \mathbf t\neq\mathbf 0.$$
Then it follows by similar arguments as for Lemma 6.1 in \cite{Pit96} that
$$P\left(\max_{\mathbf{t}\in \mathbf{g}_{u}\mathcal{K}}Z(\mathbf{t})>u+\frac{x}{u},\max_{\mathbf{t}\in \mathbf{g}_{u}( \mathfrak{\widehat{R}}_{D}\times [0,\lambda]^{m-1})}Z(\mathbf{t})>u\right)= \mathcal{H}_{D,\boldsymbol \alpha}^{x,0}(\lambda)\Psi(u)(1+o(1))$$
as $u\rightarrow\infty$, where $\mathfrak{\widehat{R}}_D=\{kD u^{-2/\alpha}: kD u^{-2/\alpha}\le \lambda, k\inn\}$ with $D>0$ is a Pickands grid in $\R$ and $\H_{D,\boldsymbol{\alpha}}^{x,0}(\lambda)$ is defined by \eqref{def: H.D}.
It also can be proved in a similar way as for Lemma 7.1 of \cite{Pit96} that
$$\H_{D,\boldsymbol \alpha}^{x,0}:=\lim_{S\rightarrow\infty}
\frac{\H_{D,\boldsymbol \alpha}^{x,0}(\lambda)}{\lambda^{m}}\in
(0,\infty).$$
\COM{and
$$\H_{D,\boldsymbol {\alpha}_0}^{0,x}=\frac{\H_{D,\alpha}^{0,x}}{\pi^{(m-1)/2}},??????$$
where $\H_{D,\alpha}^{0,x}$ is defined in Section 2.}
It is easy to check that
$$1-r(t,\mathbf{v},s,\mathbf{w})=(1+o(1))\left(|t-s|^{\alpha}+\sum_{i=1}^{m-1}\left(\frac{1}{\sqrt{2}}(w_{i}-v_{i})\right)^{2}\right)$$
as $|t-s|\rightarrow0$ and $\vnorm{\mathbf v-\mathbf w}\rightarrow 0$.
Now, by similar arguments as for Theorem 7.1 and  Corollary 7.3  of \cite{Pit96}, we have
\cL{for sufficiently large $T$ and small $b>0$}
\begin{eqnarray*}
\label{eq318}
&&\mathbb P\left\{\max_{(t,\mathbf{v})\in \mathcal D}Y(t,\mathbf{v})>v_{T}+ \frac{x}{v_{T}},
\max_{(t,\mathbf{v})\in\widetilde{\mathfrak D}}Y(t,\mathbf{v})>v_{T}\right\}\\ &&=(1+o(1))\mathbb P\left\{\max_{(t,\mathbf{v})\in [0,T^{a}]\times \mathcal S_{m-1}}Y(t,\mathbf{v})>v_{T}+ \frac{x}{v_{T}},
\max_{(t,\mathbf{v})\in\mathfrak{R}(\delta)\cap([0,T^{a}]\times\mathcal S_{m-1})}Y(t,\mathbf{v})>v_{T}\right\}\\
&&=(1+o(1)) \cdot 2^{(3-m)/2}\pi^{m/2}\ccL{(\Gamma(m/2))^{-1}}T^{a}\H_{D,\boldsymbol {\alpha}_0}^{x,0}v_{T}^{2/\alpha+m-1}\Psi(\cL{v_T}).
\end{eqnarray*}
Using further (\ref{equ413}), we get
\begin{eqnarray*}
v_{T}
&=&\frac{x+r-\sqrt{2r}\vnorm{\mathbf z}}{a_{T}}+b_{T}+o(a_{T}^{-1})\\
&=&v_{T}^{*}+b_{T}-b_{\delta, T}+(x-y)/a_{T}+o(a_{T}^{-1})\\
&=&v_{T}^{*}+\frac{\ln
\mathcal{H}_{\alpha}-\ln  \mathcal{H}_{D,\alpha} +x-y}{v_{T}^{*}}+O\left((\ln\ln (T))^{2}(\ln T)^{-3/2}\right).
\end{eqnarray*}
Observing that $v_{T}^{*}= (2\ln
T)^{1/2}(1+o(1))$, we see that the reminder $O(\cdot)$ plays a
negligible role. Therefore, using again \eqref{equ413}, 
we have
\begin{eqnarray*}
&&\mathbb P\left\{\max_{(t,\mathbf{v})\in [0,T^{a}]\times \mathcal S_{m-1}}Y(t,\mathbf{v})> v_{T},
   \max_{(t,\mathbf{v})\in\widetilde{\mathfrak D}}Y(t,\mathbf{v})> v_{T}^{*}\right\}\\
&&=2^{(3-m)/2}\pi^{m/2}\Gamma^{-1}(m/2)T^{a}\mathcal{H}_{D,\boldsymbol {\alpha}_0}^{Z_{x,y},0}(v_{T}^{*})^{2/\alpha+m-1}\Psi(v_{T}^{*})(1+o(1))\\
&&=T^{a-1}\pi^{(m-1)/2}\mathcal{H}_{D,\boldsymbol {\alpha}_0}^{Z_{x,y},0}\mathcal{H}^{-1}_{D,\alpha}e^{-y-r+\sqrt{2r}\vnorm{\mathbf z}}(1+o(1)),
\end{eqnarray*}
where $Z_{x,y}=\ln
 \mathcal{H}_{\alpha}-\ln
 \mathcal{H}_{D,\alpha} +x-y$.
Next, changing the variables in the definition of
$\mathcal{H}_{D,\boldsymbol {\alpha}_0}^{x,y}$ we get that
$\mathcal{H}_{D,\boldsymbol {\alpha}_0}^{Z_{x,y},0} \mathcal{H}_{D,\alpha} ^{-1}e^{-y}
=\mathcal{H}_{D,\boldsymbol {\alpha}_0}^{\ln
\mathcal{H}_{\alpha}+x, \ln
\mathcal{H}_{D,\alpha}+y}$, which completes the proof of the lemma. \hfill$\Box$
\\


\bigskip
{\bf Acknowledgement}: We would like to thank Enkelejd Hashorva for several valuable  suggestions and discussions. 

\bibliographystyle{plain}

 \bibliography{DiscreteChi0}

\end{document}